\font\gothic=eufm10.
\font\piccolo=cmr8.
\font\script=eusm10.
\font\sets=msbm10.
\font\stampatello=cmcsc10.
.
\def\0{{\bf 0}}
\def\1{{\bf 1}}
\def\defineq{\buildrel{def}\over{=}}
\def\EqByDef{\buildrel{\bullet}\over{=}}
\def\definiz{\buildrel{def}\over{\Longleftrightarrow}}
\def\C{\hbox{\sets C}}

\def\N{\hbox{\sets N}}
\def\R{\hbox{\sets R}}
\def\Primes{\hbox{\sets P}}
\def\T{\hbox{\sets T}}
\def\Z{\hbox{\sets Z}}
\def\square{\hbox{\vrule\vbox{\hrule\phantom{s}\hrule}\vrule}}
\def\nondiv{\!\not \; | }
\def\supporto{{\rm supp}}
\def\cf{\hbox{\stampatello c-f}\thinspace\thinspace}
\def\Gcf{G_{\hbox{\stampatello c-f}}}
\def\sf{\hbox{\stampatello s-f}\thinspace\thinspace}
\def\WA{(\hbox{\stampatello WA})}
\def\DH{(\hbox{\stampatello DH})}
\def\BH{(\hbox{\stampatello BH})}

\def\REEF{(\hbox{\stampatello R.E.E.F.})}

\def\PGrid{(\hbox{\stampatello P-Grid})}
\def\52{(5.2)}
\def\A{(\hbox{\stampatello A})}
\def\P0{(\hbox{\stampatello P}_0)}
\def\EF{(\hbox{\stampatello EF})}
\def\SingSer{\hbox{\gothic S}}
\def\Car{{\rm Car}}
\def\Win{{\rm Win}}
\def\CarT{\Car \; }
\def\WinT{\Win \; }

\def\WinPellF{\Win^{(P)}_{\ell} F\;}
\def\Pb{(P)^{\flat}}
\def\WinPellFtilde{\Win^{(P)}_{\ell} \widetilde{F}\;}
\def\Rvl{(\hbox{\stampatello Rvl})}

\def\IPP{(\hbox{\stampatello IPP})}

\def\PNT{(\hbox{\stampatello PNT})}

\def\assurdo{\rightarrow\!\leftarrow}
\def\sumflat{\mathop{{\sum}^{\flat}}}
\def\starsum{\mathop{{\sum}^{*}}}
\def\doublesum{\mathop{{\sum}{\sum}}}
\par
\centerline{\bf General elementary methods meeting elementary properties of correlations}
\bigskip
\centerline{Giovanni Coppola}\footnote{ }{MSC $2020$: $11{\rm N}05$, $11{\rm P}32$, $11{\rm N}37$ - Keywords: Ramanujan expansion, correlation, $2k-$twin primes} 
\bigskip
\bigskip
\bigskip
\par
\noindent
{\bf Abstract}. 
This is a kind of survey on properties of correlations of two very general arithmetic functions, from the point of view of {\it Ramanujan expansions}. In fact, our previous papers on these links had, as a focus, the \lq \lq {\it Ramanujan coefficients}\rq \rq \enspace of these correlations  and the resulting \lq \lq {\it R.e.e.f.}\rq \rq, i.e., {\it Ramanujan exact explicit formula}. This holds under a variety of sufficient conditions, mainly under two conditions of convergence involving correlation's {\it \lq \lq Eratosthenes Transform\rq \rq}, namely what we call {\it Delange Hypothesis} (DH) and {\it Wintner Assumption} (WA) (the former implying the latter). We proved, on Hardy-Ramanujan Journal, the \enspace {\it Hardy-Littlewood Conjecture} on {\it $2k-$twin primes}, in particular, from (DH) (implying convergence of the Ramanujan expansion, whence {\it R.e.e.f.}); then, we gave a more general proof, from (WA), again {\it entailing the R.e.e.f.}, but from another method of summation for {\it Ramanujan expansions}, that we detailed in \lq \lq On Ramanujan smooth expansions for a general arithmetic function\rq \rq. There we also started to introduce local Ramanujan coefficients. Which we deepen here, entailing a new conditional proof of quoted H-L Conjecture. 

\bigskip
\bigskip

\par
\noindent{\bf 1. Introduction. Notations, definitions \& well-known facts. Paper's Plan} 
\bigskip
\par
\noindent
As we sail through the ocean of properties of correlations, we will give for granted a series of number-theoretic very well-known properties, trying to build up, from the beginning, what we will use without quotations (but quoting equations in our calculations). We will mainly rely upon two books, for this, namely, Davenport's [D] \&  Tenenbaum's [T] books (asking help to [IKo] sometimes), where the relative proofs of quoted well-known facts may be found. In this paper, we use the term {\stampatello fact} in both senses, the common and that of \lq \lq {\stampatello small Lemma}\rq \rq. Playing with words, in case we want to briefly prove it, instead of \lq \lq Proof\rq \rq, we write \lq \lq In fact,\rq \rq. By the way, we'll use {\stampatello QED} for the end of a Proof's part, while, as usual, Proof's end will be denoted with \hfill $\square$ 
\smallskip
\par
Present \lq \lq hybrid-survey\rq \rq \thinspace deals with elementary methods linked to Ramanujan expansions; literature on the subject of correlations is a vast ocean but rather all, say, in deep sea: they're not elementary methods. And {\bf hybrid-survey means} we give an already appeared result (present Theorem 3.1 in $\S3$, see [C1]) {\bf with} a lot of {\bf  new material}, \lq \lq recent\rq \rq, meaning: not yet published. Luckily enough, we found {\bf recent} interesting results. To our results in $\S3$ we add our Th.ms 3.5, 3.6 \& Corollaries 3.2, 3.3, since present version 3. See that Corollary 3.3 proves H-L Conjecture for $2k-$twin primes when $2k$ is square-free, assuming $\A$ in $\S3$. 
\smallskip
\par
Conventional notation, we follow: {\stampatello naturals} set $\N=\{1,2,\ldots\}$, {\stampatello integers} set $\Z$, {\it primes}, esp., $p$ or $P$, filling $\Primes$; g.c.d., {\it greatest common divisor}, of $a,b\in \N$, abbreviated $(a,b)$, while $[a,b]$ is their {\it least common multiple}, l.c.m. (while, $\forall x\in \R$, the {\stampatello reals}, $[x]$ is the {\it integer part}, {\stampatello unique} $n\in \Z$: $x-1<n\le x$); whence {\it coprimality} of $a,b$, i.e. no common prime-factors, expressed as $(a,b)=1$, is codified in number theory as 
\par
\noindent {\bf Fact 1.1}
$$
\1_{(a,b)=1}=\sum_{{d|a}\atop {d|b}}\mu(d), 
\qquad 
\forall a,b\in \N
\leqno{(1.1)}
$$
\par
\noindent
giving us now opportunity, to introduce notation: 
\item{$\1$} the \lq \lq indicator function\rq \rq, like above $1_{\wp}$ {\bf indicates} $1$ in case $\wp$ is true, $0$ otherwise (abbrev., o.w.); as a consequence, it will {\bf also} indicate, in case the subscript is a set, the {\bf characteristic function of a set} of natural numbers; we agree to abbreviate with $\1:=\1_{\N}$ the characteristic function of natural numbers;
\item{$|$} the \lq \lq {\bf divides}\rq \rq, a symbol indicating, see the above, esp., $d|a$ that $a\in \Z$ is an integer multiple of $d\in \N$;
\item{$\mu$} the {\bf M\"obius function}: writing \enspace $\omega(n)\defineq |\{ p\in \Primes : p|n \}|$ \enspace for the number of $n$'s prime-factors, then $\mu(n)\defineq (-1)^{\omega(n)}$, when $n\in \Primes$ or $n=$a product of distinct primes, including $n=1$ (having 0 prime-factors) with $\mu(1)\defineq 1$, all these $n$ being called (for obvious reasons) {\it square-free numbers}; their complement in $\N$ is the set of $n$ not square-free (i.e., $\exists p\in \Primes$, $p^2|n$), on which we agree $\mu(n)\defineq 0$. Last but not least, this entails that $\mu^2(n)=|\mu(n)|$ is the {\it characteristic function} of square-free natural numbers. 
\par				
\noindent
Actually, previous Fact is a straightforward consequence of the standard [T] 
\smallskip
\par
\noindent {\bf Fact 1.2} ({\stampatello M\"{o}bius inversion})
$$
\sum_{d|n}\mu(d)=\1_{\{1\}}(n),\enspace \forall n\in \N,
\quad
\hbox{\it i.e.,}
\qquad 
\mu \ast \1 = \1_{\{1\}}\,, 
\leqno{(1.2)}
$$
\par
\noindent
giving another opportunity, to introduce our notation:\enspace $\C$ is the set of {\stampatello complex numbers}, 
\item{$\1_{\{1\}}$} is the unity in the ring of {\bf arithmetic functions}, {\stampatello a.f.}s, $f:\N\rightarrow \C$, write $f\in \C^{\N}$, with $+$ and the 
\item{$\ast$} \lq \lq {\bf Dirichlet product}\rq \rq, or Dirichlet convolution, between a.f.s [T]. 
\par
\noindent
In the following, we'll call \lq \lq {\bf M\"{o}bius inversion}\rq \rq, also, the more general statement [T]: 
$$
f=g\ast \1 \Longleftrightarrow g=f\ast \mu, 
\quad
\forall f,g:\N \rightarrow \C.
$$
\par
We recall that the $f:\N \rightarrow \C$ with $f(mn)=f(m)f(n)$, $\forall (m,n)=1$ are {\it multiplicative} (resp., $\forall m,n\in \N$, {\it completely multiplicative}); while, $f$ with $f(mn)=f(m)+f(n)$, $\forall (m,n)=1$ are {\it additive} (resp., $\forall m,n\in \N$, {\it completely additive}). Once fixed $p\in \Primes$, the {\it $p-$adic valuation} of any $n\in \N$ is $v_p(n)\defineq \max\{r\in \N_0 : p^r|n\}$ where we set $\N_0\defineq \N \cup \{0\}$. So, $\mu$ resp., $\omega$ defined above are: multiplicative resp., additive; while the {\it total number of prime-factors}, defined on all $n\in \N$ as $\Omega(n)\defineq \sum_{p|n}v_p(n)$ is completely additive, giving rise to the {\it Liouville function} $\lambda(n)\defineq (-1)^{\Omega(n)}$, a completely multiplicative a.f. The {\bf null-function}, \enspace $\0(n)\defineq 0$, $\forall n\in \N$, \enspace has all these four properties. Some Authors do not include $\0$, for example, in the multiplicative functions, but we do not follow them. The {\it von Mangoldt function} \thinspace $\Lambda(n)\defineq \log p$, if $n=p^r$ for some $r\in \N$, while $\defineq 0$ otherwise, is a kind of \lq \lq {\bf weighted characteristic function of primes}\rq \rq. For all this: [D], [T].
\smallskip
\par
Another important property we'll use is
\smallskip
\par
\noindent {\bf Fact 1.3} ({\stampatello M\"obius divisor sums into prime products})
$$
f:\N \rightarrow \C
\enspace
\hbox{\it multiplicative},
\enspace 
f\neq \0 
\quad 
\Longrightarrow 
\quad
\sum_{d|n}\mu(d)f(d)=\prod_{p|n}\left(1-f(p)\right),
\enspace
\forall n\in \N.
\leqno{(1.3)}
$$
\par
Notice that, in particular, $n=1$ gives an {\bf empty product}: always $1$ in present paper; like always $0$ the {\bf empty sum}, henceforth. Recall the definition, once fixed an $f\in \C^{\N}$, of its {\bf Eratosthenes Transform} [W]
$$
f'\defineq f\ast \mu, 
\quad 
\hbox{\it namely}
\quad
f'(d)\defineq \sum_{t|d}f(t)\mu\left({d\over t}\right),
\enspace \forall d\in \N.
$$
\par
\noindent
As a corollary to M\"obius inversion $(1.2)$, we may define it through the equivalent: $\sum_{d|n}f'(d)=f(n)$, $\forall n\in \N$. 
\medskip
\par
In order to detect a congruence condition, we may use two properties: we give the same numbering to the following two facts, distinguishing with \lq \lq a\rq \rq: for additive \& soon after with \lq \lq m\rq \rq: for multiplicative. 
\smallskip
\par
A fundamental tool is the classic 
\smallskip
\par
\noindent {\bf Fact 1.4a}({\stampatello Orthogonality of additive characters})
{\it Defining the imaginary exponential} {\stampatello modulo} $q\in \N$ 
$$
e_q(\cdot):m\in \Z \mapsto e_q(m)\defineq e^{2\pi i m/q}\in \C, 
$$
\par
\noindent
{\it we have}
$$
\1_{n\equiv 0(q)}={1\over q}\sum_{j\le q}e_q(jn), 
\enspace
\forall n\in \Z
\qquad
\Longrightarrow 
\qquad
\1_{a\equiv b(q)}={1\over q}\sum_{j\le q}e_q(j(a-b)), 
\enspace
\forall a,b\in \Z.
\leqno{(1.4a)}
$$
\par
\noindent
Here the {\bf additive character} modulo $q$, $e_q(\cdot)$, transforms sums in $\Z$ into products into the unit circle $\T\subset \C$ and $a\equiv b(q)$ abbreviates $a\equiv b(\bmod \;q)$ the {\bf congruence modulo $q$}, with the usual meaning $q|(a-b)$,here. 
\medskip
\par				
\noindent
By the way, we will call \lq \lq a $q-$periodic function\rq \rq, henceforth, a function that repeats its values after $q$ : but, we are not saying that its {\bf period}, defined to be the least natural number after which it repeats, is $q$. As an example: the $2-$periodic $\sin(\pi x)$ function has period $2$, but the $2-$periodic $\sin^2(\pi x)$ has instead period $1$. 
\medskip
\par
Another fundamental tool is the classic invention of Dirichlet, namely 
\smallskip
\par
\noindent {\bf Fact 1.4m} ({\stampatello Orthogonality of multiplicative characters: Dirichlet characters})
\par
\noindent
{\it Defining the set of $\varphi(q)$ Dirichlet characters} {\stampatello modulo} $q\in \N$ [D] {\it as the unique set of not null arithmetic functions, $q-$periodic, completely multiplicative and vanishing outside the set of naturals coprime to $q$}
$$
\chi(\cdot):n\in \N \mapsto \chi(n)\in \C, 
\quad
\hbox{\it say} \enspace \chi(\bmod \; q)
$$
\par
\noindent
{\it with its} {\bf conjugate}
$$
\overline{\chi}(\cdot):n\in \N \mapsto \overline{\chi(n)}\in \C, 
\quad
\hbox{\it say} \enspace \overline{\chi}(\bmod \; q), 
$$
\par
\noindent
{\it see }[D,Chapter 4] {\it for the explicit definition (and the characterization above, that we use as a definition), we have another way to express congruence conditions:} ({\it recall} \enspace $\Z_q^*\defineq \{n\le q : (n,q)=1\}$, {\it hereafter}) 
$$
\1_{n\equiv 1(q)}={1\over {\varphi(q)}}\sum_{\chi (\!\!\bmod q)}\chi(n), 
\enspace
\forall n\in \N
\quad
\Longrightarrow 
\quad
\1_{n\equiv a(q)}={1\over q}\sum_{\chi (\!\!\bmod q)}\overline{\chi(a)}\chi(n), 
\enspace
\forall n\in \N, \forall a\in \Z_q^*,
\leqno{(1.4m)}
$$
\par
\noindent
{\it summing over all the $\varphi(q)$ Dirichlet characters modulo $q$. Notice, however, that while additive characters detect ALL residue classes $a(\bmod \;q)$, multiplicative ones detect ONLY reduced residue classes $a\in \Z_q^*$.}
\par
\noindent
{\it For completeness, we give the other orthogonality relation, to detect the} {\bf principal character} $\chi_0(\bmod \;q)$, {\it i.e.,} $\chi_0: n\in \N \mapsto \chi_0(n)\defineq \1_{(n,q)=1}\in \C$, {\it summing over \thinspace $n(\bmod \;q)$, say } $n\le q$, {\it or equivalently over } $\Z_q^*$,
$$
\1_{\chi=\chi_0}={1\over {\varphi(q)}}\sum_{n(\!\!\bmod q)}\chi(n). 
$$
\medskip
\par
Gathering by g.c.d. and using M\"obius inversion, together with the properties of classic {\bf Euler totient} \enspace $\varphi(n)\defineq |\{m\le n : (m,n)=1\}|$, $(1.4a)$ implies $(1.8)$, see [D,Chap.26], for H\"older's \& Kluyver's formul\ae, too: (hereafter, {\it w.r.t.} abbreviates \lq \lq with respect to\rq \rq) 
\smallskip
\par
\noindent {\bf Fact 1.5} ({\stampatello Ramanujan Sums Formul\ae}) {\it Defining the} {\bf Ramanujan sum} [R] {\stampatello modulo} $q\in \N$
$$
\qquad \qquad \quad 
c_q(\cdot):m\in \N \mapsto c_q(m)\defineq \sum_{j\in \Z_q^*}\cos(2\pi jm/q)=\sum_{j\in \Z^*_q}e_q(jm)\in \C
\quad
(\hbox{\rm the \enspace \lq \lq }\sin\hbox{\rm \rq \rq \enspace function\thinspace is\thinspace odd}), 
$$
\par
\noindent
{\it recalling } $\Z_q^*$ {\it from above, we have} [D]\thinspace : $\forall a\in \Z$ {\it fixed}, $c_q(a)$ {\it is multiplicative w.r.t.} $q\in \N$ {\stampatello and the formul\ae}
\smallskip
\item{$(1.5)$} $c_q(m+q)=c_q(m)$,\enspace $\forall m\in \Z$\enspace ({\it it is $q-$periodic})\enspace {\it and }  $c_q(rn)=c_q(n)$, $\forall r\in \Z_q^*$ \enspace ({\it a variable change in } $\Z_q^*$)
\item{$(1.6)$} $c_q(n)={\displaystyle \varphi(q){{\mu(q/(q,n))}\over {\varphi(q/(q,n))}} }$, $\forall n\in \Z$ ({\stampatello H\"older formula}) [M]
\item{$(1.7)$} $c_q(n)={\displaystyle \sum_{{d|q}\atop {d|n}}d\mu(q/d) }\in \Z$, $\forall n\in \Z$ ({\stampatello Kluyver formula}) [D] 
\item{$(1.8)$} ${\displaystyle \1_{n\equiv 0(q)}=\1_{q|n}={1\over q}\sum_{d|q}c_d(n),\forall n \in \Z}$ [M]\enspace ({\it Renders} {\stampatello Ramanujan Sums} {\it ubiquitous in number theory})
\item{$(1.9)$} ([Ca], {\stampatello \lq \lq Carmichael's\rq \rq, say,) Orthogonality of Ramanujan Sums} [M] 
$$
\lim_{x\to \infty}{1\over x}\sum_{a\le x}c_q(n+a)c_{\ell}(a) = \1_{q=\ell}\cdot c_{\ell}(n), 
\qquad
\forall \ell, q\in \N, \enspace \forall n\in \Z.  
$$ 
\medskip
\par
\noindent
This Orthogonality is proved in [M] by the classic exponential sums bound, see: Theorem 1 Proof there. We'll see in $\S5$, just for fun, a Proof from Kluyver's Formula. 
\vfill
\par				
\noindent
We can come, so to speak, to the real motivation to include all these well-known properties, in particular these above for Ramanujan sums: the {\stampatello Ramanujan Expansions}. These were introduced by Ramanujan [R] soon after His definition, see previous Fact, of $c_q(a)$ : by the way, we'll call $a$ its {\it argument}; they where already known, Ramanujan main merit has been the study of His expansions, of some \lq \lq coefficients\rq \rq, say, twisted by $c_q(a)$ in a series; starting with [M] called, not by chance, {\bf Ramanujan Series} of argument $a\in \N$: 
$$
\hbox{\script R}_G(a)\defineq \sum_{q=1}^{\infty}G(q)c_q(a)
 =\lim_{x\to \infty}\sum_{q\le x}G(q)c_q(a), 
\leqno{(1.10)}
$$
\par
\noindent
assuming this limit exists in $\C$, while the {\bf coefficient} $G:q\in \N \mapsto G(q)\in \C$ is an arithmetic function giving convergence to $(1.10)$; we may also speak, for the single coefficient $G(q)$ of {\bf index} $q$, twisted by Ramanujan sum of modulus $q$, as of \lq \lq $q-${\it th Ramanujan coefficient}\rq \rq. Once we know this series converges to a complex number whatever is the argument $a\in \N$, the question is: to which function? An arithmetic function, indeed: Ramanujan started from the beginning, say, of all a.f.s, the null-function $\0$
$$
\0(a)=\sum_{q=1}^{\infty}{1\over q}c_q(a), 
\quad
\forall a\in \N
\leqno{(1.11)}
$$
\par
\noindent
that has, so to speak, a long list of consequences, we'll see in a while some, but it's interesting that $(1.11)$ is EQUIVALENT to the Prime Number Theorem, see next Fact 1.6. The main consequence is \lq \lq innocent\rq \rq, say, namely that, since $\0$ has the trivial coefficient $\0$, the null-function has {\stampatello no unique coefficients}; this innocent consequence entails a HUGE property: writing \thinspace $F=\hbox{\script R}_G$ \thinspace for \thinspace $F(a)=\hbox{\script R}_G(a)$, $\forall a\in \N$, hereafter, 
\medskip
\par
\centerline{given ANY $F:\N \rightarrow \C$, IF we have $F=\hbox{\script R}_G$, THEN $F$ has the infinitely many $F=\hbox{\script R}_{G+\lambda R_0}$}
\medskip
\par
\noindent
where $R_0(q)\defineq 1/q$ is Ramanujan's coefficient [R] of null-function in $(1.11)$ (Hardy gave $H_0(q)\defineq 1/\varphi(q)$ [H]) and $\forall \lambda \in \C$, whence an infinity of coefficients $G+\lambda R_0$ for our $F$ ! In order to deepen, compare [C1]. 
\par
Thanks to Ramanujan, we now know that {\bf each} arithmetic function has {\bf an infinity of Ramanujan expansions}, if it has one. Thanks to Hildebrand, we know that all arithmetic functions have {\stampatello at least one} Ramanujan expansion: Theorem 1.1 in [ScSp,Chapt.V]. Together, these two properties prove  that {\stampatello all arithmetic functions have each an infinity of Ramanujan Expansions.} \enspace Compare [C1], [CG]. 
\par
By the way, Schwarz \& Spilker, in their classic Book \lq \lq Arithmetical Functions\rq \rq, [ScSp], say that they take, as a definition of {\it Ramanujan coefficient}, \lq \lq the most natural\rq \rq, namely for $f:\N \rightarrow \C$ the $q-$th should be
$$
M(f\cdot c_q)/\varphi(q), 
$$
\par
\noindent
where $M(h)$ of an arithmetic function $h$ is, whenever it exists in $\C$, its {\bf mean value}, i.e. the limit
$$
M(h)\defineq \lim_{x\to \infty}{1\over x}\sum_{a\le x}h(a), 
$$ 
\par
\noindent
so that : 
$$
M(f\cdot c_q)=\lim_{x\to \infty}{1\over x}\sum_{a\le x}f(a)c_q(a). 
$$
\par
\noindent
This is a choice motivated, mostly, by the same process bringing Fourier coefficients : Orthogonality. For Ramanujan Sums, it is above $(1.9)$ ! It was discovered by Carmichael [Ca]. Compare next Definition 1.1. 
\par
The choice of Schwarz \& Spilker brings other problems, with these coefficients. In fact, take $f$ to be the characteristic function of squares: it has all $M(f\cdot c_q)=0$, so in case these would be the $f$ coefficients, in Ramanujan series, $f$ would be the null-function! See [C1,Remark 9]. A definition that is arising, from [M] : 
$$
G:\N \rightarrow \C
\enspace
\hbox{\it is\enspace a}\enspace \hbox{\stampatello Ramanujan\enspace coefficient}\enspace \hbox{\it for}
\enspace \enspace
F:\N \rightarrow \C
\quad 
\definiz
\quad
\forall a\in \N, \enspace F(a)=\hbox{\script R}_G(a)
$$
\par
\noindent
with Ramanujan series $\hbox{\script R}_G(a)$ converging $\forall a\in N$, of course. Please, don't say \lq \lq the\rq \rq: no uniqueness!  
\par				
\noindent
However, above coefficients, given thereby in [ScSp] as {\bf candidates for Ramanujan coefficients}, may be Ramanujan coefficients, with this definition. The problem is to prove it! We can't call them Ramanujan coefficients anyway (remember the case above for the characteristic function of squares), so:
\smallskip
\par
\noindent{\bf Definition 1.1} {\it We call, for a fixed arithmetic function $F$, $\forall q\in \N$, the following limit when it exists in $\C$, the $q-${\stampatello th Carmichael coefficient},indicated: }
$$
\Car_q\,F\defineq {1\over {\varphi(q)}}\lim_{x\to \infty}{1\over x}\sum_{a\le x}F(a)c_q(a),
\quad
\forall q\in \N. 
$$
\par
\noindent
The name, we chose in our previous publications, is in honor of Carmichael, for the reason we told. 
\par
Another natural candidate to be a Ramanujan coefficient was studied by Wintner [W] and the following $q-$th coefficients are a kind of \lq \lq Eratosthenian Averages\rq \rq, the name of His Book [W].
\smallskip
\par
\noindent{\bf Definition 1.2} {\it We call, for a fixed arithmetic function $F$, $\forall q\in \N$, the following series when it converges in $\C$, the $q-${\stampatello th Wintner coefficient}, indicated: }
$$
\Win_q\,F\defineq \sum_{{d=1}\atop {d\equiv 0(\!\!\bmod q)}}^{\infty}{{F'(d)}\over d},
\quad
\forall q\in \N. 
$$
\par
\noindent
This was, for Ramanujan, already a well-known way to find His coefficients (say, a heuristic). Wintner was the first to study these averages of Eratosthenes Transform : under His Assumption, following $\WA$, {\stampatello Wintner Assumption} (in His Criterion, see [ScSp])
$$
\sum_{d=1}^{\infty}{{|F'(d)|}\over d}<+\infty,
\leqno{\WA}
$$
\par
\noindent
all these coefficients exist (trivially, by comparing) AND they are exactly the above Carmichael coefficients. (As we see in [C1].)
\medskip
\par
To be more synthetic, we write:
\smallskip
\par
\noindent{\bf Definition 1.3} {\it We call, for a fixed arithmetic function $F$, in case all the $q-$th Wintner coefficients exist, the} {\stampatello Wintner Transform} {\it the arithmetic function} 
$$
\Win\;F \,:\, q\in \N \mapsto \Win_q\,F\in \C. 
$$
\par
\noindent
{\it Analogously, for a fixed arithmetic function $F$, in case all the $q-$th Carmichael coefficients exist, the} {\stampatello Carmichael Transform} {\it is the arithmetic function} 
$$
\Car\;F \,:\, q\in \N \mapsto \Car_q\,F\in \C. 
$$
\medskip
\par
Thus Wintner proved (see the above) : $\WA$ {\it for } $F$  $\Longrightarrow $ $\exists \WinT\;F$ {\it and } $\WinT\;F = \CarT\;F$. He missed (by not very far, see page's end) to prove : these, say, Carmichael/Wintner coefficients are Ramanujan coefficients. Actually, He also recognized that there are counterexamples [W], because $\WA$ is not enough to prove it. Delange, in 1976, with his Hypothesis
$$
\sum_{d=1}^{\infty}{{2^{\omega(d)}|F'(d)|}\over d}<+\infty,
\leqno{\DH}
$$
\par
\noindent
slightly stronger than $\WA$, obtained His Theorem [De1]; saying a bit more (but very important!) than the Wintner's Criterion: namely, this guarantees only  existence \& coincidence of coefficients we call, for this reason, Carmichael/Wintner, but Delange's result is: {\bf (DH)} for $F$ {\bf implies absolute convergence of Ramanujan expansion with these coefficients}. We proved $\DH$ $\Rightarrow$ Hardy-Littlewood Conjecture [C0]. Corollary 2.3 in [C2] ensures, under $\WA$, the convergence, not necessarily absolute, of Ramanujan {\bf smooth} expansion, with Carmichael/Wintner coefficients: a new summation method, here, is through \lq \lq $P-$smooth partial sums\rq \rq, then $P\to \infty$, running in the primes. We call it, of course, \lq \lq {\stampatello Wintner's Dream Theorem}\rq \rq. 
\par
Fixed $P\in \Primes$, write: $(P)\defineq \{ n\in \N : p|n \Rightarrow p\le P\}$ for $P-$smooth naturals (with $(P)$ containing $1$) and $)P(\defineq \{ n\in \N : p|n \Rightarrow p>P\}$ for $P-$sifted naturals (these contain $1$, too); both two sets appear, now, in the new third version. 
\par
We go back to history with next Fact and we will use it, in its \lq \lq M\"obius language version\rq \rq, following. 

\vfill
\eject

\par				
\noindent
We wish a break, with a moment of history of number theory to pop up, in our parade of well-known facts. 
\par
We need some basic, say, notation: recall that two functions $f(x)$ and $g(x)$ of a large real variable $x>0$ are {\it asymptotic} IFF (if \& only if), by definition, 
$$
f\sim g
\enspace \definiz \enspace
\lim_{x\to \infty}{{f(x)}\over {g(x)}}=1. 
$$
\par
It's a beautiful result, \lq \lq guessed\rq \rq, by Karl Friedrich Gauss when He was a boy; then, it was proved by Hadamard \& de la Vallee' Poussin in 1896, after decisive ideas of Riemann (1859, His \lq \lq Memoir\rq \rq, [D], with His still unproven Hypothesis); it regards prime numbers, but, apart from the unique factorization in prime-powers for natural numbers, that we learn at school, it is THE Theorem, for Prime Numbers. So, we number theorists call it PNT 
\smallskip
\par
\noindent {\bf Fact 1.6} ({\stampatello Prime Number Theorem}) {\it Defining the number of primes up to } $x\ge 1$, $x\in \R$,
$$
\pi(x)\defineq |\{p\le x : p \enspace \hbox{\stampatello prime}\}|, 
$$
\par
\noindent
{\it we have}
$$
\pi(x) \sim {x\over {\log x}},
\quad \hbox{\it as}\quad x\to \infty,
\leqno{\PNT}
$$
\par
\noindent
{\it which is equivalent } [T] {\it to the formulation, say,} {\stampatello M\"{o}bius Language PNT}
$$
\sum_{d=1}^{\infty}{{\mu(d)}\over d}=\prod_{p\in \Primes}\left(1-{1\over p}\right)=0.
\leqno{\mu\PNT}
$$
\medskip
\par
\noindent
We omit, for brevity, all the other equivalent formulations (see [T]), recalling only the two: 
$$
\sum_{n\le x}\Lambda(n) \sim \sum_{p\le x}\log p 
 \sim x,
\quad \hbox{\it as}\quad x\to \infty.
$$
\par
\noindent
By the way, in this paper we have only one {\bf logarithm}, $\log$, {\bf in} the {\bf natural base}, the Napier constant $e$ (also called {\it Euler's number}, from the beautiful {\it Euler identity} \enspace $e^{\pi i}+1=0$, with five cardinal constants). 
\smallskip
\par
We profit now to complete the notation we'll need. Landau's notation $f(x)=O(g(x))$, as $x\to \infty$ (or other variables, like $Q,N,n,$etc.) abbreviates: there is a constant, called {\bf the $O-$constant}, say $C>0$, such that for \lq \lq large\rq \rq, say, $x$ (i.e., $\forall x>x_0$, $x_0$ {\bf absolute constant}) we have : $|f(x)|\le Cg(x)$, with $g(x)>0$ trivially. The Vinogradov notation $f(x)\ll g(x)$ is perfectly equivalent to Landau's, even more practical. Both may have subscripts with the variable on which the constant depends, esp., like in the
\smallskip
\par
\centerline{
$\!\!(1.12)$\hfill 
${\displaystyle f(n)\ll_{\varepsilon} n^{\varepsilon}, 
\quad
\hbox{\it as} \enspace n\to \infty, 
 }$
\hfill ({\stampatello Ramanujan Conjecture})
}
\smallskip
\par
\noindent
equivalently expressed $f(n)=O_{\varepsilon}(n^{\varepsilon})$, indicating that the $\ll-$constant, resp., $O-$constant $C=C(\varepsilon)$ depends on $\varepsilon$. By the way, not only in this requirement above, everywhere $\varepsilon>0$ appears is both {\bf arbitrarily small} and {\bf may change from one line to the other}, especially in Proofs. Usually $x$ is a real variable, while natural numbers the others, but the difference is immaterial : our $x$ may be a large natural number, too. Last but not least, both in the $O-$ and $\ll-$ subscripts, we may have a list of variables for constant $C$. 
\par
The divisor function $d(n)\defineq {\displaystyle \doublesum_{{a\ge 1\enspace \thinspace b\ge 1}\atop {ab=n}} } 1$ satisfies Ramanujan Conjecture [T], whence: $f$ has it $\Leftrightarrow$ $f'$ has it ($f'$ is the Eratosthenes Transform of our $f$, see above). We will use this property in the following, implicitly. 
\smallskip
\par
We conclude this Analytic Number Theory briefing, with the {\stampatello Naturals' Unique Factorization} (see the above, before PNT), in its Analytic version, see [T], where, in arithmetic functions, ${\rm L} : n\in \N \mapsto \log n\in \C$
$$
\sum_{d|n}\Lambda(d)=\log n,\forall n\in \N, \hbox{\it i.e.,}\enspace \Lambda \ast \1 = {\rm L}, 
\thinspace \hbox{\it equivalent\thinspace to:} \enspace 
\Lambda'(d)=-\mu(d)\log d,\forall d\in \N, \hbox{\it i.e.,}\enspace \Lambda' = -\mu \cdot {\rm L}
\leqno{(1.13)}
$$
\par
\noindent
and, calling \lq \lq ${\rm Id}$\rq \rq \enspace the arithmetic function \lq \lq $n\in \N\mapsto n\in \C$\rq \rq, the two classic formul\ae, [T], from $(1.3)$, too: 
$$
\sum_{d|n}\varphi(d)=n,\;\forall n\in \N\enspace \hbox{\it i.e.,}\enspace \varphi \ast \1 = {\rm Id}, 
\enspace \hbox{\it equivalent\enspace to:}\enspace 
{{\varphi(n)}\over n}=\sum_{d|n}{{\mu(d)}\over d}=\prod_{p|n}\left(1-{1\over p}\right), \forall n\in \N. 
\leqno{(1.14)}
$$

\vfill
\eject

\par				
\centerline{\bf Paper's Plan}
\bigskip
\par
\noindent
After the necessary parade of classic notations, definitions \& facts, done in $\S1$, we pass at next section $\S2$, that provides the specific definitions and properties, needed in the statements of our results; that are exposed in subsequent section $\S3$; and proved in our section $\S4$, giving also further remarks and comments; last section $\S5$ contains, also, some complementary material. More in detail, 
\bigskip
\item{$\S2$} {\it Definitions and properties of general arithmetic functions and correlations}: these are necessary, for our results' statements; in particular,
\smallskip
\item\item{$\S2.1$} {\it General arithmetic functions : specific definitions \& properties}, (after recalling square-free \& cube-free properties) gives an important subset of arithmetic functions, introducing the \lq \lq $\IPP$\rq \rq, Ignoring Prime Powers a.f.s and the operation transforming $F$ into its $\IPP-$counterpart, called here, say, the \lq \lq IPPification\rq \rq, of our $F$, while
\smallskip
\item\item{$\S2.2$} {\it Correlations : initial definitions and hypotheses}, then, starts to define \& study the correlations, with a preliminary result (only one VIP, Very Important Property) that allows to SIMPLIFY ANY given, general CORRELATION; this is done, at a low price, i.e. entailing a small remainder, simply cutting the divisors for the arithmetic function that \lq \lq carries the shift\rq \rq; this allows to TRANSFORM A GENERAL \lq \lq reasonable\rq \rq (i.e., {\stampatello fair}, see there) CORRELATION INTO A BASIC HYPOTHESIS CORRELATION, having lots of good properties; then, we join, one by one, other three hypotheses, giving more \& more good features, motivated by a kind of \lq \lq imitation\rq \rq, say, of iconic \& hystory making Correlation of $2k-$twin primes, studied of course by Hardy \& Littlewood with their history-making Conjecture B [HL]; 
\medskip
\item{$\S3$} {\it Statements of results}: already appeared, for A Basic hypothesis Correlation, and most recent ones for Three-hypotheses Correlations; 
\medskip
\item{$\S4$} {\it A Brief Correlations' parade of previous result and the new ones: Proofs}, then proves these results; in particular,  
\item{$\S4.1$} {\it A Basic Hypothesis for correlations: Theorem 3.1 Proof}, then, proves Th.m $3.1$, 
\item{$\S4.2$} {\it Four Hypotheses Correlations: Proof of Theorem 3.2}, instead, gives Th.m $3.2$ Proof
\item{$\S4.3$} {\it M\"{o}bius function with Dirichlet characters: Lemmata for the Proofs of Theorems 3.3, 3.4}, then, provides necessary Lemmata and
\item{$\S4.4$} {\it Three Hypotheses results: Theorems 3.3, 3.4 Proofs}, proves the Theorems $3.3$, $3.4$, needing only the three $(0)$,$(1)$,$(2)$ hypotheses. 
\item{$\S4.5$} {\it Three Hypotheses Correlations: an assumption more gives H-L Conjecture for \sf $a=2k$}, finally proves Theorems 3.5, 3.6 \& Corollaries 3.2, 3.3. This Subsection has been added since version $3$, like next section $\S5$.
\medskip
\item{$\S5$} {\it Odds \& Ends}: contains, amongst others, some complemetary results, like the Kluyver-style Proof of: $(1.7)$ implies $(1.9)$ \& the necessary Lemma 5.1 (applied also, esp., for the Proof of Th.m 3.6); a view on future directions of research, for present arguments, concludes the paper. 

\vfill
\eject

\par				
\noindent{\bf 2. Definitions and properties of general arithmetic functions and correlations} 
\bigskip
\par
\noindent
We recall and give some definitions, with few properties, necessary to state \& prove our results at $\S3$ \& $\S4$. 
\bigskip
\par
\noindent{\bf 2.1. General arithmetic functions : specific definitions \& properties} 
\bigskip
\par
\noindent
We start recalling the definition of {\bf square-free kernel} of a natural number $a\in \N$: 
$$
\kappa(a)\defineq \prod_{p|a}p,
\quad
\forall a\in \N, 
$$
\par
\noindent
with $\kappa(1)\defineq 1$ coming from our convention, see $\S1$, that empty products are $1$. This is also called the {\stampatello radical} of our $a\in \N$. This will be used soon for the $\IPP$ arithmetic functions. In the meanwhile, recall: $q\in \N$ is {\bf cube-free} IFF, by definition, no $p$ has $p^3|q$; using $p-$adic valuation, equivalently,  
$$
q\in \N
\enspace is \enspace
\hbox{\stampatello cube-free}
\enspace \definiz \enspace
v_p(q)\le 2, 
\enspace \forall p\in \Primes. 
$$
\par
\noindent
Hereafter it will be useful to abbreviate, resp., {\stampatello square-free} with $\sf$ and {\stampatello cube-free} with {\stampatello c-f}. An easy, but decisive in formul\ae, characterization of cube-free (natural) numbers is the following: 
$$
q\in \N\enspace \cf
\thinspace \Longleftrightarrow \thinspace \thinspace 
\exists! (u,v)\in \N^2
\enspace such \enspace that: \enspace 
q=uv^2,\enspace \mu^2(u)=\mu^2(v)=1\enspace (both\enspace {\hbox{\stampatello s-f}})
\enspace and \enspace
(u,v)=1. 
\leqno{(2.1)}
$$  
\par
\noindent
In fact, $(2.1)$ is proved {\stampatello taking} (recall $p\parallel q$, $p$ divides exactly $q$, is equivalent to: $v_p(q)=1$) 
$$
u=u(q)\defineq \prod_{p\parallel q}p
\quad \hbox{\stampatello and} \quad
v=v(q)\defineq \prod_{p:v_p(q)=2}p. 
$$
\par
\noindent
Recall: a set $S\subseteq \N$ is {\bf divisor-closed} IFF, by definition, $n\in S\enspace \&\enspace d|n$ entail $d\in S$. For example, both \sf naturals and \cf ones form divisor-closed sets. 
\par
Hereafter, the superscript ${\enspace}^{\flat}$ \lq \lq {\stampatello says}\rq \rq \thinspace \sf: for summation variables, in $\sum^{\flat}$, {\stampatello and} for the elements, in $S^{\flat}$. 
\par
\noindent
In the following, we abbreviate with $\EqByDef$ , meaning \lq \lq equals by definition\rq \rq (esp. in $\S4.5$, added since version 3).
\medskip
\par
We come now to a set of arithmetic functions, with a truly arithmetic property: we call \lq \lq {\stampatello Ignoring Prime Powers}\rq \rq \thinspace (or Ignores Prime Powers), abbrev. $\IPP$, an arithmetic function $F$ that doesn't depend on $a\in \N$, but on its square-free kernel $\kappa(a)$ : 
$$
F\enspace \IPP \thinspace \definiz \thinspace \forall a\in \N, \enspace F(a)=F(\kappa(a)). 
$$ 
\par
\noindent
Of course the square-free kernel itself ignores prime powers, $\kappa$ $\IPP$, but also the number of prime factors does: $\omega$ $\IPP$. There are several examples of $\IPP$ arithmetic functions in the literature, even if this name is pretty new, as we didn't find a common \lq \lq most-used-name\rq \rq, in the literature. Maybe the most important $\IPP$ function, in number theory, is von Mangoldt's: $\Lambda$ $\IPP$. This is immediately evident from the definition, but we can use the following characterization, coming from \enspace $d|\kappa(a)$ \thinspace $\Leftrightarrow$ \thinspace $d|a$ \& $\mu^2(d)=1$ : 
$$
F\enspace \IPP \Longleftrightarrow F'=\mu^2 \cdot F'. 
$$
\par
\noindent
What does it happen when we have $F$ NOT $\IPP$, but we want another function, say $\widetilde{F}$, that agrees with $F$ on square-free values and is $\IPP$ ? Well, the \lq \lq {\stampatello IPPification}\rq \rq, of our $F$, is simply this: 
$$
\widetilde{F}(a)\defineq F(\kappa(a)), \enspace \forall a\in \N. 
$$ 
\par
\noindent
Of course previous characterization entails a very useful property, used in our Theorems Proofs, about the $\widetilde{F}$ Wintner coefficients, see $\S4.4$ : 
$$
\widetilde{F}=(\mu^2 \cdot F')\ast \1. 
$$ 
\medskip
\par
\noindent
Actually, this Ippification $\widetilde{F}$, whose Wintner coefficients we want to calculate, is not a general one: our $F$ is a {\bf correlation}, also referred to as a {\bf shifted convolution sum}. 

\vfill
\eject

\par				
\noindent{\bf 2.2. Correlations : initial definitions and hypotheses} 
\bigskip
\par
\noindent
We prefer, in our papers, to use the term, once fixed $f,g:\N \rightarrow \C$ (two arbitrary a.f.s), {\stampatello correlation} of $f$ and $g$, of {\bf length} $N\in \N$ (always fixed \& \lq \lq large\rq \rq) and {\bf shift} $a\in \N$ (also called \lq \lq argument\rq \rq) for the shifted convolution sum
$$
C_{f,g}(N,a)\defineq \sum_{n\le N}f(n)g(n+a),
\quad
\forall a\in \N,
$$
\par
\noindent
that's an a.f. depending on the shift; and $g$ is its {\bf shift-carrying factor}. We prefer \lq \lq Correlation\rq \rq, because: we are lazy (not a joke!) and we wish to avoid any link to the general theory of modular forms, where the synonym is commonly used. We recall, for this, the twice named \lq \lq elementary\rq \rq, in the title ! 
\par
By the way, our first VIP, Very Important Property, for very general Correlations is the following, \lq \lq very elementary\rq \rq, result. (We give it here, as a kind of prerequisite for our very important first Theorem, i.e. Th.m 3.1, for Correlations). It needs a very important concept, that was born in [CM], namely we may truncate the divisors of correlation's shift-carrying factor, in order to simplify, not only the correlation itself, but all its coefficients' calculations (compare this page's end for details). Hereafter, \lq \lq $\sup$\rq \rq \thinspace is usual {\it lowest upper bound}, in real analysis, while \lq \lq $\supporto$\rq \rq \thinspace will be the {\bf support} of a $h\in \C^{\N}$, i.e., $\supporto(h)\defineq \{n\in \N : h(n)\neq 0\}.$
\par
Given an arbitrary arithmetic function \enspace $g:\N \rightarrow \C$, hereafter $g\neq \0$, in case \enspace $\sup \;\supporto(g')\in \N$, namely the Eratosthenes Transform of our $g$ has a finite support, we define \enspace $g$ a {\bf truncated divisor sum}, with top divisor defined as: $\max \supporto(g')$. Any integer $Q$ that is greater than or equal to the top divisor will be called a {\bf range}, for $g$; while the top divisor itself will be called the {\bf exact range} of our $g$. We will not need, however, the concept of exact range and, also, it is not immediate, of course, to calculate (for the truncated divisor sums $g$). Next result allows to pass, with a small error, from shift-carrying $g$ very general (any arithmetic function!) to a $Q-$truncated couterpart $g_Q$ that is, instead, a {\stampatello truncated divisor sum} ! 
\par
Once fixed $Q\in \N$ and given any $g:\N \rightarrow \C$, we define {\bf its $Q-$truncated} divisor sum, simply cutting all the divisors $d>Q$ :
$$
g_Q(n)\defineq \sum_{{d|n}\atop {d\le Q}}g'(d), 
\quad 
\forall n\in \N. 
$$
\par
\noindent
This is a truncated divisor sum with range $Q$ (and exact range $\le Q$), that agrees with $g$ over all $n\le Q$. 
\smallskip
This property briefly proves the following result, compare [C0] and [C1]. 
\smallskip
\par
\noindent{\bf Theorem 2.1}. ({\stampatello The divisors' Cut}). 
\par
\noindent
{\it Let $f,g:\N \rightarrow \C$ have } $\Vert f\Vert_{\infty}\defineq \max_{n\le N}|f(n)|$,\enspace $\Vert g'\Vert_{\infty}\defineq \max_{N<d\le N+a}|g'(d)|$. {\it Then, uniformly } $\forall a\in \N$, 
$$
C_{f,g}(N,a)=C_{f,g_N}(N,a)+O(\Vert f\Vert_{\infty}\thinspace \Vert g'\Vert_{\infty}\thinspace a). 
$$
\par
\noindent
{\it Furthermore, if } $f$ \& $g$ {\it satisfy} {\stampatello Ramanujan Conjecture} $(1.12)$, {\it too}, {\it this entails} ({\it again, uniformly} $\forall a\in \N$)
$$
C_{f,g}(N,a)=C_{f,g_N}(N,a)+O_{\varepsilon}(N^{\varepsilon}(N+a)^{\varepsilon}a).
$$
\smallskip
\par
\noindent {\bf Proof}. We prove more general formula in one line: 
$$
C_{f,g}(N,a)-C_{f,g_N}(N,a)=\sum_{N<q\le N+a}g'(q)\sum_{{n\le N}\atop {n\equiv -a(\!\!\bmod q)}}f(n)
\ll \Vert f\Vert_{\infty}\thinspace \Vert g'\Vert_{\infty}\sum_{N<q\le N+a}1
\ll \Vert f\Vert_{\infty}\thinspace \Vert g'\Vert_{\infty}\thinspace a. 
\enspace 
\enspace 
\enspace 
\square
$$
\medskip
\par
\noindent
This result allows to assume: $\exists Q\defineq \max(\supporto(g'))$ and $Q\le N$. Notice that we may also have an arithmetic function $g(m)$, say,  mimicking the divisor function $d(m)$ through the \lq \lq flipping of divisors\rq \rq, say $g(m):=2\sum_{d|m\, ,\, d\le \sqrt{N}}1$; this is ALREADY a truncated divisor sum: $Q=[\sqrt N]$. However, the importance of above Theorem is the universality w.r.t. $g$, namely it holds FOR ALL $g:\N \rightarrow \C$ and uniformly for all $f:\N \rightarrow \C$. 
\medskip
\par
Joining a technical requirement, say, our correlation has to be {\stampatello fair} , we get the \lq \lq Basic Hypothesis\rq \rq, for our correlation, following. See [C1] for details about technical problems in calculating Carmichael coefficients for general correlations, NOT satisfying next $\BH$. Further details in $\S4.1$, Remark 1. 
\medskip
\par				
\noindent
We give now the definition of {\stampatello Basic Hypothesis} for Correlations, abbrev. \lq \lq $\BH$\rq \rq, that we will use often. 
\medskip
\par
\centerline{\bf Basic Hypothesis:}
\smallskip
\par
\noindent
The {\stampatello correlation} of $f:\N \rightarrow \C$ and $g_Q:\N \rightarrow \C$, defined as (recall from above), once fixed {\stampatello length} $N\in \N$, 
\item{} ${\displaystyle C_{f,g_Q}(N,a)\defineq \sum_{n\le N}f(n)g_Q(n+a), \forall a\in \N, }$ as a function of {\stampatello shift} $a\in \N$, \enspace {\stampatello has two requirements}: 
\item{$1.$} The function $g_Q$ is {\stampatello of} $\underline{\hbox{\stampatello range}}$ \enspace $Q\in \N$, $Q\le N$, i.e. $\supporto(g'_Q)\subset [1,Q]$ : \enspace ${\displaystyle g_Q(m)=\sum_{{d|m}\atop {d\le Q}}g'(d), }$ \enspace where the function $g:\N \rightarrow \C$ has $g'=g'_Q$ {\stampatello or its Eratosthenes Transform $g'$ is restricted to $g'_Q$, vanishing outside} $[1,Q]$. (In fact, we saw in $\S2$ {\stampatello the divisors' cut} how to do it.)
\item{$2.$} Then, the {\stampatello correlation} $C_{f,g_Q}(N,a)$ {\stampatello is } $\underline{\hbox{\stampatello fair}}$ : meaning by definition that {\stampatello the $a-$dependence} ({\stampatello shift-dependence of correlation}) {\stampatello is } $\underline{\hbox{\stampatello only}}$ {\stampatello in the } $g_Q$ {\stampatello argument}, see the above definition. We are {\stampatello excluding other dependencies}, of our $f$ or of our $g_Q$, or of one or both of their supports. 
\par 
\noindent 
We show a list of subsequent, consecutive Hypotheses for our \enspace $C_{f,g_Q}$\enspace and we will provide three of them : $(1)$,$(2)$,$(3)$, see the following; for this reason, Basic Hypothesis, abbr. $\BH$, will be abbreviated as well as hypothesis $(0)$: it is a kind of granted, by the divisors' cut Theorem 2.1. We shortly give $(1)$,$(2)$,$(3)$. 
\item{$(1)$} $g_Q$ $\IPP$
\item{$(2)$} $f=\1_{\Primes}\cdot f$
\item{$(3)$} $f$ and $g_Q$ {\stampatello satisfy Ramanujan Conjecture}
\par
We give brief comments: $g_Q$ $\IPP$ describes an arithmetic property of our shift carrying factor, but has no immediate effect on arithmetic properties of \enspace $C_{f,g_Q}$, with respect to the shift; saying \enspace $f=\1_{\Primes}\cdot f$, of course, means: $f$ is {\stampatello supported on primes}; finally, we may assume $(3)$ for $f$ and $g_Q$, since we may think to be able to {\stampatello renormalize both factors} of our Correlation. 
\par 
Further remarks \& comments, $\S4.2$ beginning, give {\stampatello deep reasons to introduce these Hypotheses}. 

\bigskip

\par
\noindent{\bf 3. Statements of results} 
\bigskip
\par
\noindent
We state our first result, for $\BH-$correlations. This is the only one, already appeared elsewhere (see [C1]). 
\smallskip
\par
\noindent {\bf Theorem 3.1}. ({\stampatello All \lq \lq Basic Hypothesis\rq \rq$-$Correlations' main properties})
\par
\noindent{\it Let } $C_{f,g_Q}(N,a)$ {\it be a $\BH-$correlation, i.e., satisfying } 1. {\it and } 2. {\it in the above} {\stampatello Basic Hypothesis}. {\it Then,} 
\smallskip
\item{A.} {\it It is periodic, whence bounded, in all natural shifts } $a\in \N$, 
\item{B.} {\it whence, it has both Carmichael coefficients and the coincident Wintner coefficients, both vanishing on indices beyond } $Q$. 
\item{C.} {\it For these coefficients, we have } {\stampatello the explicit formula}, {\it in terms of $g_Q$ Ramanujan coefficients $\widehat{g_Q}(\ell)$}: 
$$
\Win_{\ell}\enspace C_{f,g_Q}(N,\bullet) = \Car_{\ell}\enspace C_{f,g_Q}(N,\bullet) 
 = {{\widehat{g_Q}(\ell)}\over {\varphi(\ell)}}\; \sum_{n\le N}f(n)c_{\ell}(n),
\qquad
\forall \ell \in \N.  
$$
\par
\noindent
A $\BH$ $C_{f,g_Q}(N,a)$ has {\stampatello Analytic Part} a {\it finite Ramanujan Expansion} with above coefficients, defined as: 
$$
A_{f,g_Q}(N,a)\defineq \sum_{q\le Q}{{\widehat{g_Q}(q)}\over {\varphi(q)}}\; \sum_{n\le N}f(n)c_q(n)c_q(a),
\enspace 
\forall a\in \N. 
\enspace 
\left(\hbox{\stampatello For}\enspace \widehat{g_Q}(\ell)\enspace \hbox{\stampatello definition:\thinspace Th.m\thinspace 3.1\thinspace Proof}\right)
$$ 
\par
\noindent
Sometimes we call \lq \lq Basic Correlation\rq \rq, hereafter, a $\BH-$Correlation with Ramanujan Conjecture for $f,g_Q$.
\par
We pass to our second result, proved by a \lq \lq tour\rq \rq, say, of our Four Hypotheses above in $\S2.2$. 
\smallskip
\par
\noindent {\bf Theorem 3.2}. ({\stampatello \lq \lq Four Hypotheses Theorem\rq \rq \enspace About Basic Correlations})
\par
\noindent
{\it Let } $f:\N\rightarrow\C$, $g_Q:\N\rightarrow\C$ {\it and their } {\stampatello correlation} $C_{f,g_Q}(N,a)$ {\it satisfy the } {\stampatello Four Hypotheses}:
\item{$(0)$} {\it The } {\stampatello correlation} $C_{f,g_Q}(N,a)$ {\it fulfills } {\stampatello Basic Hypothesis}
\item{$(1)$} {\it The } {\stampatello Ramanujan coefficient \enspace $\widehat{g_Q}$ \enspace is square-free supported} : $\widehat{g_Q}=\mu^2 \cdot \widehat{g_Q}$; {\it equivalently, we assume } $g_Q$ \IPP
\item{$(2)$} {\it The } {\stampatello function $f$ is supported on primes}: $f=\1_{\Primes} \cdot f$
\item{$(3)$} {\stampatello Both} {\it the functions} $f$ {\stampatello and} $g_Q$ {\stampatello satisfy Ramanujan Conjecture} $(1.12)$
\smallskip
\par				
\noindent
{\it Then} 
$$
\forall a\in \N, \enspace C_{f,g_Q}(N,a)={\hbox{\script P}}_{f,g_Q}(N,a)+{\hbox{\script S}}_{f,g_Q}(N,a), 
\enspace \hbox{\it say,} \enspace \hbox{\lq \lq {\stampatello primary\enspace and\enspace secondary\enspace parts}\rq \rq,}
$$
\par
\noindent
{\it where}
\par
\leftline{${\displaystyle 
{\hbox{\script P}}_{f,g_Q}(N,a)\defineq \sum_{q\le Q}\mu(q)\widehat{g_Q}(q)\sum_{m|q}\mu(m)m\sum_{{p\le N}\atop {{(p,m)=1}\atop {-p\equiv a(\!\!\bmod m)}}}f(p)
 }$
}
\par
\rightline{${\displaystyle 
=\sum_{q\le Q}\mu(q)\widehat{g_Q}(q)\sum_{m|q}{{\mu(m)m}\over {\varphi(m)}}\sum_{\chi\atop {(\!\!\bmod m)}}\sum_{p\le N}f(p)\overline{\chi(-p)}\chi(a), 
\quad
\forall a\in \N 
 }$
}
\par
\noindent
{\it and}
$$
{\hbox{\script S}}_{f,g_Q}(N,a)\defineq \sum_{p|a}f(p)p\sum_{q'}\mu(q')\widehat{g_Q}(pq')\sum_{{m'|q'}\atop {{a\over p}\equiv -1(\!\!\bmod m')}}\mu(m')m'
 = O_{\varepsilon}\left(N^{\varepsilon}a\right), 
\quad
\forall a\in \N.
$$ 
\medskip
\par
\noindent
Notice in the above statement of Hypothesis $(1)$, in previous Theorem, the equivalence
$$
g_Q \IPP \enspace \Longleftrightarrow \enspace \widehat{g_Q}=\mu^2 \cdot \widehat{g_Q} 
$$
\par
\noindent
that will be proved, then, in the beginning of $\S4.2$, well before the Proof of Theorem 3.2. 
\medskip
\par
The interest in this \lq \lq new object\rq \rq, the IPP-ification of a fixed arithmetic function (while the literature is plenty of $\IPP$ arithmetic functions, even if they don't yet have an official name!), comes from the fact that, for example, the Hardy-Littlewood Conjecture (see [C0] for a conditional Proof) has an $\IPP$ singular series.
\par
By the way, we have no idea, if someone has already UNCONDITIONALLY proved that, say, Hardy-Littlewood correlation (the one with shift $a=2k$, for $2k-$twin primes) is \lq \lq close\rq \rq, say, to its IPP-ification. Of course we have lots of data \& Heuristics about that, indicating that this correlation is \lq \lq close to be $\IPP$\rq \rq, i.e. \enspace $C_{\Lambda,\Lambda}(N,a)\sim \widetilde{C_{\Lambda,\Lambda}}(N,a)$, \enspace at least, for fixed $a=2k$. We are building a general theory for $\IPP$ a.f.s and for IPP-ification of not $\IPP$ ones; still, we don't know how to prove, and under which conditions, that \enspace $C_{f,g_Q}(N,a)\sim \widetilde{C_{f,g_Q}}(N,a)$. We trust the possibility to give an elementary, unconditional Proof.
\medskip
\par
Our third (Th.m 3.3) and, resp., fourth (Th.m 3.4) result gives, actually, explicit formul\ae, for both the principal part and, resp., the secondary part for the IPPification of a, say, Three Hypotheses Correlation, namely one having $(0)$,$(1)$,$(2)$.  
\medskip
\par
\noindent
We start with the first coefficients' calculation, that for Primary Part. 
\smallskip
\par
\noindent {\bf Theorem 3.3}. ({\stampatello Wintner's coefficients for Primary Part Ippification})
\par
\noindent
{\it Let } $f:\N\rightarrow\C$, $g_Q:\N\rightarrow\C$ {\it and their } {\stampatello correlation} $C_{f,g_Q}(N,a)$ {\it satisfy first Three Hypotheses} $(0)$, $(1)$, $(2)$ {\it in previous Theorem} 3.2. {\it This result defines the } {\stampatello Primary Part} \enspace ${\hbox{\script P}}_{f,g_Q}(N,a)$. {\it Then} 
$$
\Win_{\ell}\; \widetilde{{\hbox{\script P}}_{f,g_Q}} = {{\widehat{g_Q}(\ell)}\over {\varphi(\ell)}}\sum_{p\le N}f(p)c_{\ell}(p)
 +{{\mu(\ell)\widehat{g_Q}(\ell)}\over {\varphi(\ell)}}\sum_{p|\ell}f(p)\varphi(p)
  -{{\mu(\ell)}\over {\varphi(\ell)}}\sum_{p\nondiv \ell}f(p)\widehat{g_Q}(p\ell), 
\quad
\forall \ell \in \N. 
$$
\medskip
\par
We give now the coefficients for the Secondary Part. Next formula RHS (Right Hand Side) is \lq \lq small\rq \rq. 
\smallskip
\par
\noindent {\bf Theorem 3.4}. ({\stampatello Wintner's coefficients for Secondary Part Ippification})
\par
\noindent
{\it Let } $f:\N\rightarrow\C$, $g_Q:\N\rightarrow\C$ {\it and their } {\stampatello correlation} $C_{f,g_Q}(N,a)$ {\it satisfy first Three Hypotheses} $(0)$, $(1)$, $(2)$ {\it in previous Theorem} 3.2. {\it This result defines the } {\stampatello Secondary Part} \enspace ${\hbox{\script S}}_{f,g_Q}(N,a)$. {\it Then} 
$$
\Win_{\ell}\; \widetilde{{\hbox{\script S}}_{f,g_Q}} = -{{\mu(\ell)\widehat{g_Q}(\ell)}\over {\varphi(\ell)}}\sum_{p|\ell}f(p)\varphi(p)
 +{{\mu(\ell)}\over {\varphi(\ell)}}\sum_{p\nondiv \ell}f(p)\widehat{g_Q}(p\ell), 
\quad
\forall \ell \in \N. 
$$
\medskip
\par
\noindent
We wish a bit to surprise ourselves: adding together previous two Theorems RHS we get the surprise, that correlation's Wintner coefficients are THE SAME AS Wintner's coefficients for correlation's IPP-ification !!! 
\smallskip
\par
\noindent {\bf Corollary 3.1}. ({\stampatello Wintner's coefficients of Three-Hypotheses-Correlations Ippification})
\par
\noindent
{\it Under the first Three Hypotheses $(0),(1),(2)$ of Theorem 3.2 for our Correlation \enspace $C_{f,g_Q}(N,a)$, we obtain the coincidence of Wintner Transform of } $C_{f,g_Q}(N,a)$ {\it and  Wintner Transform of its  IPPification } $\widetilde{C_{f,g_Q}}(N,a)$. 

\vfill
\eject

\par				
\noindent
We give now our most recent results, proved in $\S4.5$ (a new subsection), for present $3$rd version of our paper. 
\medskip
\par
The first new result we present is still for a $3-$hypotheses Correlation $C_{f,g_Q}$, having IPPification \enspace $\widetilde{C_{f,g_Q}}$; it regards the local \lq \lq $P-$Wintner Transforms\rq \rq: the $P-$smooth Wintner coefficients, whence these transforms, are defined in $\S4.5$; say: same Hypotheses $(0),(1),(2)$ entail \enspace $\lim_P\; \Win^{(P)} \widetilde{C_{f,g_Q}}=\WinT \widetilde{C_{f,g_Q}}$ (we express this more precisely, see next Theorem). This is not very surprising, as its Proof comes \lq \lq adapting\rq \rq, so-to-speak, Lemma 4.2 arguments used for classic Wintner coefficients, to these $P-$local versions : this is done, very easily, in Lemma 4.3 (compare this, in $\S4.5$, and Lemma 4.2, both in next $\S4$). Precise statement follows. 
\medskip
\par
Once given Theorem 3.3 and Theorem 3.4, following their Proofs easily entails the following result. 
\smallskip
\par
\noindent {\bf Theorem 3.5}. ({\stampatello $P-$Wintner's coefficients for Three-Hypotheses-Correlations Ippification})
\par
\noindent
{\it Let $f$, $g_Q$ and their Correlation\enspace $C_{f,g_Q}(N,a)$\enspace satisfy the first Three Hypotheses $(0),(1),(2)$ of Theorem 3.2. This result defines the, resp.,}  {\stampatello primary, ${\hbox{\script P}}_{f,g_Q}(N,a)$, and secondary, ${\hbox{\script S}}_{f,g_Q}(N,a)$, part} {\it of our Correlation. Then}
$$
\lim_P \Win^{(P)}_{\ell} \widetilde{{\hbox{\script P}}_{f,g_Q}} = \Win_{\ell} \; \widetilde{{\hbox{\script P}}_{f,g_Q}}
\quad \hbox{\stampatello and} \quad
\lim_P \Win^{(P)}_{\ell} \widetilde{{\hbox{\script S}}_{f,g_Q}} = \Win_{\ell} \; \widetilde{{\hbox{\script S}}_{f,g_Q}},
\enspace \hbox{\stampatello both} \enspace \hbox{\stampatello uniformly} \enspace
\forall \ell \in \N. 
$$
\par
\noindent
{\stampatello By linearity and}\quad ${\displaystyle \widetilde{C_{f,g_Q}}(N,a)=\widetilde{{\hbox{\script P}}_{f,g_Q}}(N,a)+\widetilde{{\hbox{\script S}}_{f,g_Q}}(N,a) }$, \enspace $\forall a\in \N$,\quad {\stampatello these two entail}
$$
\lim_P \Win^{(P)}_{\ell} \widetilde{C_{f,g_Q}}=\Win_{\ell} \; \widetilde{C_{f,g_Q}},
\quad \hbox{\stampatello uniformly}  \enspace \forall \ell \in \N. 
\leqno{(3.1)}
$$
\par
\noindent
The reason why $(3.1)$ above is so important may be appreciated only \lq \lq renormalizing\rq \rq, i.e. considering $\widetilde{F}(a):=\widetilde{C_{f,g_Q}}(N,a)/N$ ($\forall a,N\in \N$) : assuming this \lq \lq renormalization\rq \rq, so-to-speak, has a Ramanujan Expansion (restricted to square-free $a\in \N$) with Ramanujan Coefficient $G:\N \rightarrow \C$ converging absolutely, $(3.1)$ LHS (Left Hand Side) may be calculated in terms of $G$; roughly speaking, this is the content, applied to $F$ being the renormalized $3-$Hypotheses Correlation as a particular case, of our next result: it calculates $\lim_P \Win^{(P)}_{\ell} \widetilde{F}$ in terms of $G$. \hfill Recall: $P-${\stampatello smooth} naturals, $(P)$, and $P-${\stampatello sifted} ones, $)P($, defined at p.5. 
\medskip
\par
Under an assumption on Ramanujan Expansion's $G$ for $\widetilde{F}$, an Explicit Formula links $\Win^{(P)} \widetilde{F}$ to $G$. 
\smallskip
\par
\noindent {\bf Theorem 3.6}. ({\stampatello $P-$Wintner's coefficients for $\widetilde{F}$ in terms of its Ramanujan Coefficients})
\par
\noindent
{\it Let } $F:\N \rightarrow \C$ {\it have } {\stampatello IPPification $\widetilde{F}$, with the following property:} ({\it regarding $\widetilde{F}$ Ramanujan Expansion, converging absolutely on \sf arguments} $a\in \N$, {\it recall}\enspace $\N^{\flat}$ {\it are \sf naturals} \& \cf {\it cube-free ones}) 
$$
\exists G:\N \rightarrow \C
\enspace \hbox{\it such that }
\enspace \underline{\hbox{\stampatello both}} \enspace
\forall a\in \N^{\flat},\enspace \widetilde{F}(a)=\sum_{q=1}^{\infty}G(q)c_q(a)
\enspace \underline{\hbox{\stampatello and}} \enspace
\sum_{q\; \cf}|G(q)|<\infty. 
$$
\par
\noindent
{\it Then, } $\forall P\in \Primes$ {\stampatello fixed } {\it and } $\forall \ell \in \Pb$, {\stampatello i.e., $\forall \ell \in (P)$ square-free}, 
$$
\Win_{\ell}^{(P)}\; \widetilde{F}=\sum_{s|\ell}\sum_{{t\in (P)}\atop {(t,\ell)=1}}\mu(s)\mu(t)\sum_{r\in )P(}\mu(r)G(\ell rst^2).
\leqno{(3.2)}
$$
\par
\noindent
({\stampatello Notice: in above RHS, $\ell rst^2$ is a cube-free natural number}.)
{\it In turn, } {\stampatello above formula entails}
$$
\lim_P \thinspace \Win_{\ell}^{(P)}\; \widetilde{F}=\sum_{s|\ell}\sum_{{t=1}\atop {(t,\ell)=1}}^{\infty}\mu(s)\mu(t)G(\ell st^2), 
\quad \hbox{\stampatello uniformly} \enspace \forall \ell \in \N^{\flat}.
$$
\par
\noindent
Before to apply it, we {\stampatello summarize three results above for $3-$hyp.$^{\hbox{{\stampatello ses}}}$ Correlations} $C_{f,g_Q}(N,a)$ : 
\par
From Corollary 3.1, we know : \enspace $\WinT \widetilde{C_{f,g_Q}}=\WinT C_{f,g_Q}$
\par
From Theorem 3.1, the {\stampatello formula}: \enspace $\Win_{\ell}\; C_{f,g_Q}={{\widehat{g_Q}(\ell)}\over {\varphi(\ell)}}\sum_{p\le N}f(p)c_{\ell}(p)$, $\forall \ell \in \N$
\par 
From Theorem 3.5, we also get : \enspace $\lim_P \thinspace \Win_{\ell}^{(P)} \widetilde{C_{f,g_Q}}=\Win_{\ell}\; \widetilde{C_{f,g_Q}}$, {\stampatello uniformly} $\forall \ell \in \N$
\medskip
\par
\noindent
In all, a Three-hp.s Correlation $C_{f,g_Q}(N,a)$ has:\enspace ${\displaystyle \lim_P \thinspace \Win_{\ell}^{(P)} \widetilde{C_{f,g_Q}}={{\widehat{g_Q}(\ell)}\over {\varphi(\ell)}}\sum_{p\le N}f(p)c_{\ell}(p) }$, {\stampatello unif.}$^{\hbox{{\stampatello ly}}}$ \thinspace $\forall \ell \in \N$. 

\vfill
\eject

\par				
\noindent
We {\stampatello need following Assumption}: 
$$
\exists G : \N \rightarrow \C
\enspace \hbox{\it such\thinspace that} \enspace 
\forall a\in \N \enspace \sf \hbox{\stampatello fixed},
\enspace \widetilde{C_{f,g_Q}}(N,a)=N\sum_{q=1}^{\infty}G(q)c_q(a) 
\enspace \hbox{\stampatello and} \enspace \sum_{q\; \cf}|G(q)|<\infty. 
\leqno{\A}
$$
\par
\noindent
We will prove Hardy-Littlewood Conjecture for $a=2k$ \sf, assuming the \lq \lq Absolute Convergence\rq \rq, for the Ramanujan Coefficients (of {\it normalized Correlation}), in $\A$, for a suitable 3-hp.s Correlation having $Q=N$.
\medskip
\par
However, we want to assume as less as possible, to get interesting general results for 3-hp.s Correlations. In such\enspace $C_{f,g_Q}$, we need to  define {\stampatello Wintner's Prime} (compare [C2] for a general definition), i.e., write, say, $P_0=P_0(C_{f,g_Q})\defineq \min\{ P\in \Primes : \supporto(\WinT C_{f,g_Q})\subseteq (P)\}$\enspace which depends, apart from $f,g$, on $Q$ \& $N$: we will {\stampatello assume} 
$$
P_0\to \infty,
\quad \hbox{\stampatello as} \enspace
N\to \infty.
\leqno{\P0}
$$
\medskip
\par
Under Assumption $\A$, all 3-hp.s Correlations $C_{f,g_Q}(N,a)$ with $\P0$ above have an Explicit Formula linking $G$ in $\A$ with $\widehat{g_Q}$ and $f$. 
\smallskip
\par
\noindent {\bf Corollary 3.2}. ({\stampatello Functions $f$, $\widehat{g_Q}$ for 3-hp.s $\widetilde{C_{f,g_Q}}$ in terms of its Ramanujan Coefficients})
\par
\noindent
{\it Let } $C_{f,g_Q}(N,a)$ {\it be } {\stampatello a three-hypotheses Correlation satisfying both assumptions $\A$ \& $\P0$.} {\it Then,} {\stampatello taking}\enspace $F(a):=C_{f,g_Q}(N,a)/N$, {\it we have trivially} $\widetilde{F(a)}\EqByDef \widetilde{C_{f,g_Q}}(N,a)/N$, $\forall a,N\in \N$
{\it which satisfies the property required in Theorem 3.6, since it comes from $\A$. Hence} 
$$
{1\over N}\lim_P \thinspace \Win_{\ell}^{(P)}\; \widetilde{C_{f,g_Q}}=\lim_P \thinspace \Win_{\ell}^{(P)}\; \widetilde{F}=\sum_{s|\ell}\sum_{{t=1}\atop {(t,\ell)=1}}^{\infty}\mu(s)\mu(t)G(\ell st^2), 
\quad \hbox{\stampatello uniformly} \enspace \forall \ell \in \N^{\flat}.
$$
\par
\noindent
{\it Using the three results, for 3-hp.s Correlations,} {\stampatello quoted above},
$$
\sum_{s|\ell}\sum_{{t=1}\atop {(t,\ell)=1}}^{\infty}\mu(s)\mu(t)G(\ell st^2)={{\widehat{g_Q}(\ell)}\over {N\varphi(\ell)}}\sum_{p\le N}f(p)c_{\ell}(p), 
\enspace \forall \ell \in \N\enspace \hbox{\stampatello square-free}. 
\leqno{\EF}
$$
\par
\noindent
{\it This } {\stampatello Explicit Formula, applying again $\A$ \& $\P0$, entails:} ({\it recall \enspace $A_{f,g_Q}(N,a)$ \enspace def.at page $10$})
$$
\forall a\in N\enspace \sf \hbox{\stampatello fixed}, 
\quad
\widetilde{C_{f,g_Q}}(N,a)=C_{f,g_Q}(N,a)=A_{f,g_Q}(N,a)+o(N), 
\enspace \hbox{\stampatello as } P_0\to \infty, \enspace \hbox{\stampatello when } N\to \infty.
$$
\medskip
\par
\noindent
Corollary 3.2 is a consequence of 4 results: one Corollary (i.e., 3.1) and three Theorems (i.e., 3.1, 3.5, 3.6). 
\medskip
\par
\noindent
In turn, next result is a kind of immediate consequence of Corollary 3.2. 
\medskip
\par
We can prove H-L Conjecture for square-free shifts $a=2k\in \N$, only assuming $\A$. 
\smallskip
\par
\noindent {\bf Corollary 3.3}. ({\stampatello Hardy-Littlewood Conjecture for $a=2k$ square-free follows assuming $\A$})
\par
\noindent
{\it Let } {\stampatello the three-hypotheses Correlation} $C_{f,g_Q}(N,a)=C_{\mu^2\,\Lambda,\Lambda_N}(N,a)$ {\it satisfy above assumption $\A$.} {\it Then,}
$$
\forall a\in \N\enspace \sf \hbox{\stampatello fixed}, 
\quad
\widetilde{C_{\mu^2\,\Lambda,\Lambda_N}}(N,a)=C_{\mu^2\,\Lambda,\Lambda_N}(N,a)=A_{\mu^2\,\Lambda,\Lambda_N}(N,a)+o(N), 
\enspace \hbox{\stampatello as } N\to \infty.
$$
\par
\noindent
{\it This } {\stampatello entails Hardy-Littlewood Conjecture for \sf $a=2k-$twin primes:}
$$
\forall a\in \N\enspace \sf \hbox{\stampatello fixed}, 
\quad
C_{\Lambda,\Lambda}(N,a)=\hbox{\gothic S}(a)N+o(N), 
\enspace \hbox{\stampatello as } N\to \infty.
$$
\medskip
\par
\noindent
We recall that \enspace $\hbox{\gothic S}$ \enspace is the {\stampatello classic singular series}, that we define as (and write in [C0]: $\hbox{\gothic S}_{\Lambda,\Lambda}$)
$$
\hbox{\gothic S}(a)\defineq \sum_{q=1}^{\infty}{{\mu^2(q)}\over {\varphi^2(q)}}c_q(a),
\quad
\forall n\in \N,
$$
\par
\noindent
whence we have the classic Euler product, proving, in particular, that if $a\in \N$ is odd, then \enspace $\hbox{\gothic S}(a)=0$.

\vfill
\eject

\par				
\noindent{\bf 4. A Brief Correlations' parade of previous result and the new ones: Proofs} 
\bigskip
\par
\noindent
\par
\noindent{\bf 4.1 A Basic Hypothesis for correlations: Theorem 3.1 Proof} 
\bigskip
\par
\noindent
We dedicate this subsection to our most fundamental hypothesis, which will be even completed, adding other hypotheses; for example, see next subsections  on the Three and Four Hypotheses Correlations. 
\par
\noindent
\par
\noindent {\bf Proof of Th.m 3.1}. By above definition of $g_Q$, we obtain {\stampatello its Ramanujan expansion, with Wintner coefficients}:
$$
\forall m\in \N, \thinspace g_Q(m)\buildrel{(1.8)}\over{=\!=\!=}\sum_{d\le Q}{{g'(d)}\over d}\sum_{\ell | d}c_{\ell}(m)
 =\sum_{\ell \le Q}\widehat{g_Q}(\ell)c_{\ell}(m),
\quad
\hbox{\stampatello with} \enspace 
\widehat{g_Q}(\ell)\defineq \sum_{{d\le Q}\atop {d\equiv 0(\!\!\bmod \ell)}}{{g'(d)}\over d}
$$ 
\par
\noindent
and notice that, by definition, these {\stampatello coefficients vanish outside} $[1,Q]$ : i.e., $\supporto(\widehat{g_Q})\subset [1,Q]$. 
\par
This formula and correlation's definition entail : 
$$
C_{f,g_Q}(N,a) = \sum_{q\le Q}\widehat{g_Q}(q)\sum_{n\le N}f(n)c_q(n+a), 
\enspace 
\forall a\in \N
\leqno{(\ast)}
$$
\par
\noindent
whence, {\stampatello from $q-$periodicity of} $c_q$, see $(1.5)$, {\stampatello the {\script Q}$-$periodicity of correlation follows} (compare [C3] for finer aspects of a $\BH$ Correlation's periodicity): 
$$
C_{f,g_Q}(N,a+\hbox{\script Q}) = C_{f,g_Q}(N,a), 
\enspace 
\forall a\in \N
$$
\par
\noindent
where \enspace {\script Q}$\defineq $l.c.m.$(2,3,\ldots,Q)$ : as a function of $a\in \N$, the $\BH-$correlation $C_{f,g_Q}(N,a)$ is $\hbox{\script Q}-\underline{\hbox{\stampatello periodic}}$ with {\stampatello period} {\tt T}$\in \N$, {\stampatello say}, where {\tt T}$|\hbox{\script Q}$. In particular, a $\BH-$correlation is {\bf bounded} (on all $a\in \N$). Hence, 
$$
C'_{f,g_Q}(N,d)\defineq \sum_{t|d}C_{f,g_Q}(N,t)\mu\left({d\over t}\right), 
\quad 
\forall d\in \N
$$
\par
\noindent
(by definition), its Eratosthenes Transform {\stampatello has Ramanujan Conjecture} $(1.12)$ : 
$$
C'_{f,g_Q}(N,d)\ll_{N,Q} \sum_{t|d}1\ll_{N,Q,\varepsilon} d^{\varepsilon}, 
\quad 
\forall d\in \N
\leqno{(4.1)}
$$
\par
\noindent
because : $C_{f,g_Q}(N,a)\ll_{N,Q} 1, \forall a\in \N$ and the divisor function [T], too, has Ramanujan Conjecture $(1.12)$ [T]. See that $(4.1)$ follows the property given soon after $d(n)$ definition, page 6. 
\par
Being bounded, we may apply the Theorem and Remark 1.5 of [De2] to get, in particular (compare Proposition 3 in [C1]):
\smallskip 
\par
\centerline{IF the Carmichael transform of a $\BH-$correlation $C_{f,g_Q}(N,a)$ exists,} 
\par
\centerline{THEN there exists its Wintner transform, too, and they are equal.}  
\smallskip
\par
As a consequence, once we know Carmichael coefficients of a $\BH-$correlation, we also have for free its Wintner coefficients. These are, say, uncomfortable to calculate from the definition, due to the general difficulty of correlation's Eratosthenes Transform. 
\par
We only need, after these considerations, to add {\stampatello Carmichael's Orthogonality} $(1.9)$ to \hfill conclude:
\par
\noindent
Definition $1.1$, of Carmichael's $\ell-$th coefficient, using $(\ast)$ above and Orthogonality $(1.9)$ gives, $\forall \ell \in \N$ fixed, 
$$
\Car_{\ell}\enspace C_{f,g_Q}(N,\bullet) = {1\over {\varphi(\ell)}}\lim_{x\to \infty}{1\over x}\sum_{a\le x}C_{f,g_Q}(N,a)c_{\ell}(a)
 \buildrel{(\ast)}\over{=\!=\!=}{1\over {\varphi(\ell)}}\widehat{g_Q}(q)\sum_{n\le N}f(n)\lim_{x\to \infty}{1\over x}\sum_{a\le x}c_q(n+a)c_{\ell}(a)
$$
$$
\buildrel{(1.9)}\over{=\!=\!=}{1\over {\varphi(\ell)}}\widehat{g_Q}(\ell)\sum_{n\le N}f(n)c_{\ell}(n).
$$
\par
\noindent 
This proves the point C. explicit formula above; whence, point B.'s vanishing, from:\enspace $\supporto(\widehat{g_Q})\subset [1,Q]$. \hfill $\square$ 
\medskip
\par				
\noindent {\bf Remark 1.} See that the hypothesis (that seems hidden) of having a fair correlation is in $(\ast)$ above, during the Proof, and it's the main reason why $\BH-$correlations may benefit from Carmichael's Orthogonality. Compare correlations above, say, before the cut (namely, truncation of divisors). Not only the range of $g$ may be $+\infty$, so to speak (as $a\to \infty$) : we have problems in exchanging sums, because the $g-$Ramanujan coefficients indices may depend on $a$, no fairness, at all! So to speak, we \lq \lq designed\rq \rq, for these reasons, the two requirements in $\BH$ above.\hfill $\diamond$
\medskip
Even before the introduction of $\BH$, we gave [CM,Theorem 1] the (Carmichael) coefficients formula in Theorem 3.1, applying \lq \lq Carmichael's Formula\rq \rq, see Definition 1.1, under the two hypotheses of \lq \lq purity\rq \rq, of Ramanujan expansion, and its uniform convergence, see [CM]. Actually, compare quoted Theorem 1, we already gave some characterizations for the {\stampatello Ramanujan exact explicit formula}, the $\REEF$, to hold. This is the {\it finite Ramanujan expansion} : 
$$
C_{f,g_Q}(N,a)=\sum_{\ell \le Q}\left({{\widehat{g_Q}(\ell)}\over {\varphi(\ell)}}\sum_{n\le N}f(n)c_{\ell}(n)\right)c_{\ell}(a),
\quad
\forall a\in \N,  
\leqno{\REEF}
$$
\par
\noindent 
which is, actually, a {\bf fixed length} expansion, as it is finite, with a number of terms NOT DEPENDING ON $a\in \N$. Theorem 1 [CM] in its hypotheses has the finite {\it range} for $g$ (but we didn't define it). Also, the concept of {\it purity} of Ramanujan expansion in [CM] has evolved into the concept of {\it fair} correlation. \hfill See [C1]. 
\medskip
\par
We see the $\REEF$'s RHS to be the {\stampatello analytic part} of our $C_{f,g_Q}(N,a)$, defined at page 10: 
$$
A_{f,g_Q}(N,a)\EqByDef \sum_{\ell \le Q}\left({{\widehat{g_Q}(\ell)}\over {\varphi(\ell)}}\sum_{n\le N}f(n)c_{\ell}(n)\right)c_{\ell}(a),
\quad
\forall a\in \N,  
$$
\par
\noindent 
where the {\stampatello coefficients are Carmichael-Wintner's} for the Correlation. This definition above may be generalized to all $F:\N \rightarrow \C$ with a {\stampatello finitely supported Wintner Tansform}, as we do in [C2].  
\par
We use\enspace $A_{f,g_Q}(N,a)$\enspace above in Corollaries 3.2 \& 3.3, added since present version 3. (See that, from our Corollary 3.1, $A_{f,g_Q}=\widetilde{A_{f,g_Q}}$, for the {\stampatello three-hypotheses correlations}.)   
\medskip
\par
\noindent 
The new idea applied in [C2], then, is to introduce the {\stampatello Ramanujan smooth expansion}: $\forall a\in \N$, 
$$
\hbox{\script R}^{\piccolo smooth}_G(a)\defineq \lim_{P\to \infty}\sum_{q\in (P)}G(q)c_q(a), 
\enspace 
\hbox{\it limit\enspace running\enspace over} \enspace P\in \Primes, 
\leqno{(4.2)}
$$ 
\par
\noindent 
where the single {\bf smooth partial sums} are over the set $(P)$ of $P-${\it smooth naturals} (defined at page 5) and {\bf are finite}, from the (recall $v_p$, the $p-$adic valuation, defined at page 2)
$$
c_q(a)\neq 0 \enspace \Longrightarrow \enspace v_p(q)\le v_p(a) + 1, \forall p|q.
\leqno{({\hbox{\stampatello Ramanujan vertical limit}})}
$$
\par
\noindent 
({\stampatello Actually}, \lq \lq $\Longleftrightarrow$\rq \rq \thinspace holds here, as we have proved in [C2].) We'll abbreviate $\Rvl$, for above property. 
\medskip
\par
\noindent
In general, $(4.2)$ is very different from $\hbox{\script R}_G(a)$ (esp., $\hbox{\script R}^{\piccolo smooth}_{\1}=\0$, but $\hbox{\script R}_{\1}$ doesn't converge), but in case $G$ is finitely supported, of course $\hbox{\script R}_G=\hbox{\script R}^{\piccolo smooth}_G$; whence, choosing $G=\Win\enspace C_{f,g_Q}$, finitely supported by above Theorem 3.1, Corollary 2.3 [C2] gives the above $\REEF$, under $\WA$ for the correlation; whence under this $\WA$ we may prove H-L Conjecture, compare [C0]. 
\par
Actually, we got in [C0] the $\REEF$ for $C_{\Lambda,\Lambda_N}(N,a)$, whence H-L Conjecture for $a=2k-$twins, {\stampatello under $\DH$}, see $\S1$ above, Delange Hypothesis (stronger than $\WA$, in fact). 

\vfill
\eject

\par				
\noindent{\bf 4.2. Four Hypotheses Correlations: Proof of Theorem 3.2} 
\bigskip
\par
\noindent
Hereafter we prove our recent results.\hfill Hypotheses $(0)$,$(1)$,$(2)$,$(3)$ named here are defined in $\S2.2$. 
\par
\noindent
We will now specialize our study of $\BH-$correlations, joining three hypothesis, to the Basic Hypothesis, that we count as $(0)$, say for granted. We are going to join, in order of importance, the Hypothesis $(1)$ assuming that the Ramanujan Coefficients $\widehat{g_Q}(q)$ are square-free supported. This comes both as a great simplification in our handling of $q$, mainly from $(d,q/d)=1$ for all $d|q$, and as a conceptual property of $g_Q$. In fact, we defined in $\S2.1$ the arithmetic functions $\IPP$, Ignoring Prime Powers, as those $F$ with $F(a)$ depending not on $a$, but on its square-free kernel; these $F$ are exactly those with square-free supported Eratosthenes Transform $F'$ : compare $\S2.1$ characterization. 
\par
\noindent
Going back to our $g_Q$, then, we may see that $g'_Q = \mu^2 \cdot g'_Q$ IFF $\widehat{g_Q} = \mu^2 \cdot \widehat{g_Q}$; in fact, if we call {\it divisors} $d$ the arguments of $g'_Q$ and {\it moduli} $q$ those of $\widehat{g_Q}$, when all $d$ are square-free it follows that, by $\widehat{g_Q}$ definition, also all $q$ are square-free: vice-versa, the finite Ramanujan expansion of $g_Q$ in $\S4.1$ allows to calculate $g'_Q$ as: 
$$
g'_Q(d)=d\sum_{K\le Q/d}\mu(K)\widehat{g_Q}(dK),
\quad
\forall d\in \N, 
$$
\par
\noindent
from $(1.7)$ and $(1.2)$, and this tells us that all $q=dK$ square-free, here, entail all $d$ square-free, too. In all, 
$$
g_Q\, \IPP
\quad \Longleftrightarrow \quad
g'_Q = \mu^2 \cdot g'_Q
\quad \Longleftrightarrow \quad
\widehat{g_Q} = \mu^2 \cdot \widehat{g_Q}. 
$$
\par
\noindent
Thus $(1)$ has two forms, depending on opportunity: we may require $g_Q\, \IPP$, for theoretical aspects, OR $\widehat{g_Q} = \mu^2 \cdot \widehat{g_Q}$, when making calculations. For this reason, we gave the equivalence of these two forms (proved now, proving previous double equivalence) in Theorem 3.2 statement. 
\par
The theoretical aspects come from: assuming that the {\bf shift-carrying factor} in the correlation, i.e., $g_Q$ in $C_{f,g_Q}(N,a)$ (it carries the shift $a\in \N$ since it appears inside $g_Q(n+a)$, by definition) {\bf doesn't depend on prime powers} entails that, if our correlation has a heuristic formula which is a kind of Hardy-Littlewood main term for $C_{f,g_Q}(N,a)$, say: (see the classic [HL]) 
$$
C_{f,g_Q}(N,a)=N\SingSer_{f,g_Q}(a)+\hbox{\it good\enspace error}
$$
\par
\noindent
then the $f,g_Q-${\bf singular series} $\SingSer_{f,g_Q}(a)$ $\IPP$ (like for the Classic Singular Series of Hardy-Littlewood). This, in turn, will furnish $\SingSer_{f,g_Q}(a)$ a Ramanujan expansion with square-free coefficients (they're exactly $\mu^2(q)/\varphi^2(q)$, for quoted Classic). Compare [CM] for a deeper explanation, involving \lq \lq singular sums\rq \rq. 
\par
In fact, in the Hardy-Littlewood case, $g_Q=\Lambda_N$, after cutting divisors at $Q=N$; the original $\Lambda$ is $\IPP$ (being, if not $0$ since $n$ is a power of $p$, $\log p$ on ALL $p$ PRIME POWERS...) and $\Lambda_N$ is a truncated divisor sum, from it, whence $\IPP$, too; another way to see it, is looking at Eratosthenes Transform: if square-free supported at the beginning, once we restrict it, in this case truncating its support, after it remains square-free supported ! 
\par
In some sense, we are trying to get closer and closer, with our hypotheses, to Hardy-Littlewood case.
\par
In next Theorems' calculations, it will be clear that, for the $f$ support, all coprimality conditions may be simplified; the most \lq \lq drastic\rq \rq, say, restriction on $f$ support, then, is to embed it into primes ! Thus, Hypothesis $(2)$ is in fact $f$ prime supported. Needless to say, while $g_Q$ is $\IPP$, now $f$ is vanishing outside primes so, we confirm it, we are even closer to Hardy-Littlewood case. (But, for example, $f$ may miss a lot of primes, or, on the contrary, it may weight all primes too much, and so on.) 
\par
Last but not least, at some point (the latest, the better), we have to make bounds of smaller or less important, say, error terms. Then, Ramanujan Conjecture, our Hypothesis $(3)$, making our two correlation factors $f$ and $g_Q$ be not too big, is a kind of must. Also this is matched by Hardy-Littlewood case, that has even lighter factors: bounded by $\log N$ each (instead of $\ll_{\varepsilon} N^{\varepsilon}$). 
\medskip
\par
However these Four Hypotheses, in the following applied rather everywhere in this subsection, are a compromise; between generality, as we're not too sticked to Hardy-Littlewood case, and simplification, especially of calculations. 
\smallskip
\par
Recall ($\S3$) that we speak informally of \lq \lq Basic Correlations\rq \rq, meaning with $\BH$ and hypothesis $(3)$. 

\vfill
\eject

\par				
\noindent {\bf Proof of Th.m 3.2}. We give a \lq \lq tour\rq \rq, say, of our Four Hypotheses, see $\S2.2$, and use Kluyver's $(1.7)$, with the property : $p|m$ \& $p\equiv -a(\bmod m)$ $\Rightarrow $ $p|a$ and Dirichlet's characters Orthogonality $(1.4m)$: 
$$
C_{f,g_Q}(N,a)\buildrel{(0)}\over{=\!=\!=}\sum_{q\le Q}\widehat{g_Q}(q)\sum_{n\le N}f(n)c_q(n+a)
 \buildrel{(1)}\over{=\!=\!=}\sum_{q\le Q}\mu(q)\widehat{g_Q}(q)\sum_{m|q}\mu(m)m\sum_{{n\le N}\atop {-n\equiv a(\!\!\bmod m)}}f(n)
$$
$$
\buildrel{(2)}\over{=\!=\!=}\sum_{q\le Q}\mu(q)\widehat{g_Q}(q)\left(\sum_{m|q}\mu(m)m\sum_{{p\le N}\atop {{(p,m)=1}\atop {-p\equiv a(\!\!\bmod m)}}}f(p)+\sum_{{p|a}\atop {p|q}}f(p)\mu(p)p\sum_{m'\left|{q\over p}\right.\,,\,m'\left|\left({a\over p}+1\right)\right.}\mu(m')m'\right)
$$
$$
\buildrel{(3)}\over{=\!=\!=}{\hbox{\script P}}_{f,g_Q}(N,a)+O_{\varepsilon}\left(Q^{\varepsilon}a\sum_{q\le Q}{1\over q}\sum_{p|q}\sum_{m'\left|{q\over p}\right.\,,\,(m',p)=1}1\right)
$$
$$
\buildrel{(0)}\over{=\!=\!=}\sum_{q\le Q}\mu(q)\widehat{g_Q}(q)\sum_{m|q}{{\mu(m)m}\over {\varphi(m)}}\sum_{\chi\atop {(\!\!\bmod m)}}\sum_{p\le N}f(p)\overline{\chi(-p)}\chi(a)+O_{\varepsilon}\left(N^{\varepsilon}a\sum_{q\le Q}{{d(q)}\over q}\right),
$$
\par
\noindent
whence the formula above for {\stampatello primary part} and above definition and bound for {\stampatello secondary part}: 
$$
\sum_{q\le Q}{{d(q)}\over q}=\sum_{m\le Q}\sum_{n\le Q/m}{1\over {mn}}
 \ll (\log Q)\sum_{m\le Q}{1\over m}
  \ll \log^2 Q 
   \ll_{\varepsilon} Q^{\varepsilon}, 
$$
\par
\noindent
where ${\displaystyle d(q)=\doublesum_{{m,n\ge 1}\atop {mn=q}}1 }$ is the {\stampatello divisor function}, see $\S1$, then use again $Q\le N$ from $(0)$.\hfill $\square$ 
\medskip
\par
\noindent {\bf Remark 2.} We call \lq \lq {\stampatello primary}\rq \rq \enspace and \enspace \lq \lq {\stampatello secondary}\rq \rq \enspace these parts, not only to say that \thinspace ${\hbox{\script S}}_{f,g_Q}(N,a)$ \thinspace may be \lq \lq {\stampatello somehow neglected}\rq \rq, by bound above; mainly, harmonics $(\bmod \enspace m)$ are {\stampatello primary}, while $(\bmod \enspace m')$ are {\stampatello secondary}.\hfill $\diamond$

\medskip

\par
\noindent
The Theorems $3.3$ \& $3.4$ simply calculate Wintner coefficients for the, say, Ippifications of Primary \& Secondary Parts, defined and started to study, in previous result. 

\medskip

In order to prove these two Theorems, we need next two Lemmata. 

\vfill
\eject

\par				
\noindent{\bf 4.3. M\"{o}bius function with Dirichlet characters: Lemmata for the Proofs of Theorems 3.3, 3.4} 
\bigskip
\par
\noindent 
We need a Lemma that gives the principal character contribution, to our Theorems $3.3$, $3.4$ calculations. 
\par
We call it the Pinch Lemma, for two reasons: formal, it \lq \lq pinches\rq \rq, somehow, the two variables $d$ and $m$, entailing (in next Theorems' Proofs) $\ell | q$, whence: $\ell \le Q$, a precious property; informal, it's so beautiful, also for consequences, that I had to pinch myself to convince I really found it!
\smallskip
\par
\noindent {\bf Lemma 4.1} ({\stampatello Pinch Lemma}) {\it Let } $d,m\in \N$. {\it Then, for } $\chi_0(\bmod m)$, 
$$
\sum_{t|d}\mu(t)\chi_0(t)=\sum_{{t|d}\atop {(t,m)=1}}\mu(t)
 =\1_{\kappa(d)|m}. 
$$
\smallskip
\par
\noindent {\bf Proof}. We give three arguments. First one uses $(1.3)$ finite Euler products : ({\stampatello empty prod.s}$=1$) 
$$
\sum_{{t|d}\atop {(t,m)=1}}\mu(t)=\prod_{{p|d}\atop {p\nondiv m}}\left(1-1\right)
 =\1_{p|d\,\Rightarrow\,p|m}
  =\1_{\kappa(d)|m}.
$$
\hfill QED
\par
\noindent
Second, write, from $(1.1)$, $\1_{(t,m)=1}={\displaystyle \sum_{{n|t}\atop {n|m}}\mu(n) }$, whence, by $(1.2)$, we end up again with $\1_{\kappa(d)|m}$ : 
$$
\sum_{{t|d}\atop {(t,m)=1}}\mu(t)=\sum_{{t|\kappa(d)}\atop {(t,m)=1}}\mu(t)
 =\sum_{{n|\kappa(d)}\atop {n|m}}\mu(n)\sum_{{t|\kappa(d)}\atop {t\equiv 0(\!\!\bmod n)}}\mu(t)
  =\sum_{{n|\kappa(d)}\atop {n|m}}\mu(n)\sum_{\widetilde{t}\left|{{\kappa(d)}\over n}\right.}\mu(n)\mu(\widetilde{t})
   =\sum_{{n|\kappa(d)}\atop {n|m}}\mu^2(n)\sum_{\widetilde{t}\left|{{\kappa(d)}\over n}\right.}\mu(\widetilde{t}). 
$$
\hfill QED
\par
\noindent
Third, set \enspace $n=(\kappa(d),m)$, so
$$
\sum_{{t|d}\atop {(t,m)=1}}\mu(t)=\sum_{{t|\kappa(d)}\atop {(t,m)=1}}\mu(t)
 =\sum_{t\left|{\kappa(d)\over n}\right.}\mu(t)
  =\1_{n=\kappa(d)}
   =\1_{(\kappa(d),m)=\kappa(d)}
    =\1_{\kappa(d)|m}
$$
\par
\noindent
again by M\"obius inversion $(1.2)$.\hfill $\square$ 
\medskip
\par
\noindent
We need another Lemma for the non-principal characters contribution, to our Theorems $3.3$, $3.4$ calculations. 
\par
It blends together PNT in M\"{o}bius language and following bound for non-principal Dirichlet characters: 
$$
\sum_{p\le x}{{\chi(p)}\over p}\ll_m 1,
\enspace \hbox{\rm as} \enspace x\to \infty,
\quad 
\forall \chi\neq \chi_0(\bmod \;m), 
\leqno{(4.3)}
$$
\par
\noindent
a form of Dirichlet's Theorem in Arithmetic Progressions (proved by Mertens: [D],Ch.7; compare [IKo],$\S2.3$). 
\medskip
\par
\noindent {\bf Lemma 4.2}. ({\stampatello PNT in M\"{o}bius language for non-principal Dirichlet characters})
\par
\noindent
{\it Let } $p\in \Primes$, $\ell\in \N$. {\it Fix } $m\in \N$ {\it and fix any non-principal Dirichlet character } $\chi\neq \chi_0(\bmod \;m)$. {\it Define}
$$
\Sigma(\chi,\ell,m)\defineq \sum_{d\equiv 0(\!\!\bmod \ell)}{{\mu(d)}\over d}\sum_{t|d}\mu(t)\chi(t), 
$$
\par
\noindent
{\it while, fixed } $m'\in \N$ {\it and fixed any non-principal Dirichlet character } $\chi\neq \chi_0(\bmod \;m')$, {\it define}
$$
\Sigma_p(\chi,\ell,m')\defineq \sum_{{d\equiv 0(\!\!\bmod \ell)}\atop {d\equiv 0(\!\!\bmod p)}}{{\mu(d)}\over d}\sum_{t'\left|{d\over p}\right.}\mu(t')\chi(t'). 
$$
\par
\noindent
{\it Then}
$$
\Sigma(\chi,\ell,m)=\Sigma_p(\chi,\ell,m')=0. 
$$
\smallskip
\par				
\noindent {\bf Proof}. We start with $\Sigma(\chi,\ell,m)$, setting $d=\ell \widetilde{d}$, to get, by multiplicative properties of $\mu,\chi$ and their product's convolution with $\1$, still multiplicative, using $(1.3)$ twice, 
$$
\Sigma(\chi,\ell,m) = {{\mu(\ell)}\over \ell}\sum_{t|\ell}\mu(t)\chi(t)\sum_{(\widetilde{d},\ell)=1}{{\mu(\widetilde{d})}\over {\widetilde{d}}}\sum_{t|\widetilde{d}}\mu(t)\chi(t)
 = {{\mu(\ell)}\over \ell}\sum_{t|\ell}\mu(t)\chi(t)\prod_{\widetilde{p}\nondiv \ell}\left( 1-{{1-\chi(\widetilde{p})}\over {\widetilde{p}}}\right).
$$
\par
\noindent
This Euler product is, factoring as follows, 
$$
\prod_{\widetilde{p}\nondiv \ell}\left( 1-{{1-\chi(\widetilde{p})}\over {\widetilde{p}}}\right)=\prod_{\widetilde{p}|\ell}\left( 1-{1\over {\widetilde{p}}}\right)^{-1}\cdot \prod_{\widetilde{p}\in \Primes}\left( 1-{1\over {\widetilde{p}}}\right)\cdot \prod_{\widetilde{p}\nondiv \ell}\left( 1+{{\chi(\widetilde{p})}\over {\widetilde{p}-1}}\right)=0, 
$$
\par
\noindent
since for $\chi\neq \chi_0(\bmod \;m)$, we have (we exclude, now, the $\widetilde{p}=2-$factor, vanishing if $\chi(2)=-1$ and $\ell$ is odd, making the product with $\,\chi\,$ vanish!) 
$$
\prod_{\widetilde{p}>2,\widetilde{p}\nondiv \ell}\left( 1+{{\chi(\widetilde{p})}\over {\widetilde{p}-1}}\right)=\exp\left(\sum_{\widetilde{p}>2,\widetilde{p}\nondiv \ell}\log \left( 1+{{\chi(\widetilde{p})}\over {\widetilde{p}-1}}\right)\right)\ll \exp\left(\sum_{\widetilde{p}>2,\widetilde{p}\nondiv \ell}{{\chi(\widetilde{p})}\over {\widetilde{p}}}\right)\ll_m 1, 
$$
\par
\noindent
from $(4.3)$, while the other factors' product vanishes, from $\mu\PNT$.\hfill QED
\par
\noindent
For $\Sigma_p(\chi,\ell,m')$, the situation is similar, but we have to distinguish two cases: $p|\ell$ aut $p\not |\ell$ (\lq \lq aut\rq \rq, here, is the {\bf exclusive or} of course). We start with $p|\ell$, giving
$$
\Sigma_p(\chi,\ell,m')=\sum_{d\equiv 0(\!\!\bmod \ell)}{{\mu(d)}\over d}\sum_{t'\left|{d\over p}\right.}\mu(t')\chi(t')
 ={{\mu(\ell)}\over {\ell}}\sum_{(\widetilde{d},\ell)=1}{{\mu(\widetilde{d})}\over {\widetilde{d}}}\sum_{t'\left|\left({\widetilde{d}\cdot {\ell \over p}}\right)\right.}\mu(t')\chi(t')
$$
$$
={{\mu(\ell)}\over {\ell}}\sum_{t'\left|{\ell \over p}\right.}\mu(t')\chi(t')\sum_{(\widetilde{d},\ell)=1}{{\mu(\widetilde{d})}\over {\widetilde{d}}}\sum_{t'\left|\widetilde{d}\right.}\mu(t')\chi(t')
={{\mu(\ell)}\over {\ell}}\sum_{t'\left|{\ell \over p}\right.}\mu(t')\chi(t')\prod_{\widetilde{p}\nondiv \ell}\left( 1-{{1-\chi(\widetilde{p})}\over {\widetilde{p}}}\right)
 =0, 
$$
\par
\noindent
since now again $\widetilde{d}=d/\ell$, which is coprime to $\ell/p$ here, and we use vanishing product above.
\par
\noindent
In case $p$ doesn't divide $\ell$, set $d=p\ell \widetilde{d}$, getting
$$
\Sigma_p(\chi,\ell,m')=\sum_{d\equiv 0(\!\!\bmod p\ell)}{{\mu(d)}\over d}\sum_{t'\left|{d\over p}\right.}\mu(t')\chi(t')
 ={{\mu(p\ell)}\over {p\ell}}\sum_{t'|\ell}\mu(t')\chi(t')\sum_{(\widetilde{d},p\ell)=1}{{\mu(\widetilde{d})}\over {\widetilde{d}}}\sum_{t'|\widetilde{d}}\mu(t')\chi(t')
$$
$$
={{\mu(p\ell)}\over {p\ell}}\sum_{t'|\ell}\mu(t')\chi(t')\prod_{\widetilde{p}\nondiv \ell,\widetilde{p}\neq p}\left( 1-{{1-\chi(\widetilde{p})}\over {\widetilde{p}}}\right)
 =0, 
$$
\par
\noindent
because: this Euler product may differ from previous at most on one factor, whence vanishing again.\hfill $\square$
\medskip
\par
We want to apply these two Lemmata to the calculation of Wintner coefficients, not of our correlation $C_{f,g_Q}(N,a)$, that we already did in Theorem 3.1; but, of its \lq \lq {\stampatello IPPification}\rq \rq, that is, the arithmetic function (of shift $a\in \N$) that is the same of $C_{f,g_Q}(N,a)$  on square-free shifts $a\in \N$, being an $\IPP-$arithmetic function, too. A moment's reflection reveals that (like, in general, at $\S2.1$) we have to take our function and multiply its Eratosthenes Transform by $\mu^2$ (say, \lq \lq killing\rq \rq \enspace all non-square-free divisors), then we sum back over divisors: in our case, 
$$
\widetilde{C_{f,g_Q}}(N,a)\defineq \sum_{d|a}\mu^2(d)C_{f,g_Q}'(N,d),
\qquad 
\forall a\in \N. 
$$
\par
Why are we so interested, in a kind of \lq \lq obscure object\rq \rq? Well, by definition it agrees with the original function, our $C_{f,g_Q}(N,a)$ here, on square-free numbers. Second, it is $\IPP$ and this suggests something, for our correlations: we are assuming that they satisfy, for the Four Hypotheses, in particular $(1)$: $g_Q$ $\IPP$; this, from $\BH$ (Hypothesis $(0)$, recall) and Theorem 3.1, gives: the Carmichael/Wintner coefficients of our correlation are square-free supported ! In other words, since by definition also the Wintner coefficients of our correlation's Ippification are square-free supported, we may reasonably expect that they are intimately linked to correlation's coefficients ! Time to calculate the Wintner coefficients of $\widetilde{C_{f,g_Q}}$ and next two Theorems will do it: the $\widetilde{\enspace}$ operator is linear and we have above Theorem 3.2 decomposition; that doesn't need $(3)$. The Theorems $3.3$, $3.4$ and their Corollary $3.1$, see $\S3$, will skip $(3)$, as they all don't have error terms. 

\vfill
\eject

\par				
\noindent{\bf 4.4. Three Hypotheses results: Theorems 3.3, 3.4 Proofs}
\bigskip
\par
\noindent {\bf Proof of Th.m 3.3}. We start from Theorem 3.2 formula for ${\hbox{\script P}}_{f,g_Q}(N,a)$, once fixed $a\in \N$; by definition of Eratosthenes Transform $\S1$, $\forall d\in \N$, 
$$
{\hbox{\script P}}_{f,g_Q}'(N,d)=\sum_{t|d}\mu\left({d\over t}\right)\sum_{q\le Q}\mu(q)\widehat{g_Q}(q)\sum_{m|q}{{\mu(m)m}\over {\varphi(m)}}\sum_{\chi\atop {(\!\!\bmod m)}}\sum_{p\le N}f(p)\overline{\chi(-p)}\chi(t)
$$
\par
\noindent
whence, using $\mu(d/t)=\mu(d)\mu(t)$ for all square-free $d$ and any $t|d$, 
$$
\mu^2(d){\hbox{\script P}}_{f,g_Q}'(N,d)=\mu(d)\sum_{q\le Q}\mu(q)\widehat{g_Q}(q)\sum_{m|q}{{\mu(m)m}\over {\varphi(m)}}\sum_{\chi\atop {(\!\!\bmod m)}}\sum_{p\le N}f(p)\overline{\chi(-p)}\sum_{t|d}\mu(t)\chi(t)
$$
\par
\noindent
and, recalling $\S2.1$ formula and the Definition 1.2 of Wintner $\ell-$th coefficient, abbreviated as \enspace $\Win_{\ell}\; \widetilde{{\hbox{\script P}}_{f,g_Q}}$ , 
$$
\Win_{\ell}\; \widetilde{{\hbox{\script P}}_{f,g_Q}} = \sum_{d\equiv 0(\!\! \bmod \ell)}{{\mu(d)}\over d}\sum_{q\le Q}\mu(q)\widehat{g_Q}(q)\sum_{m|q}{{\mu(m)m}\over {\varphi(m)}}\sum_{\chi\atop {(\!\!\bmod m)}}\sum_{p\le N}f(p)\overline{\chi(-p)}\sum_{t|d}\mu(t)\chi(t) 
$$
$$
=\sum_{q\le Q}\mu(q)\widehat{g_Q}(q)\sum_{m|q}{{\mu(m)m}\over {\varphi(m)}}\sum_{\chi\atop {(\!\!\bmod m)}}\sum_{p\le N}f(p)\overline{\chi(-p)}\sum_{d\equiv 0(\!\! \bmod \ell)}{{\mu(d)}\over d}\sum_{t|d}\mu(t)\chi(t), 
$$
\par
\noindent
which, by Lemma 4.2, has no terms with non-principal characters, because the $d-$sum here is \enspace $\Sigma(\chi,\ell,m)$ \enspace there defined; with  the terms with principal characters a non-vanishing contribute, from Lemma 4.1, give  
$$
\Win_{\ell}\; \widetilde{{\hbox{\script P}}_{f,g_Q}} = \sum_{q\le Q}\mu(q)\widehat{g_Q}(q)\sum_{m|q}{{\mu(m)m}\over {\varphi(m)}}\sum_{p\le N,(p,m)=1}f(p)\sum_{{d\equiv 0(\!\! \bmod \ell)}\atop {d|m}}{{\mu(d)}\over d}
$$
\par
\noindent
which, by 
$$
\sum_{{d\equiv 0(\!\! \bmod \ell)}\atop {d|m}} {{\mu(d)}\over d} = \1_{\ell|m} {{\mu(\ell)}\over \ell}\cdot 
 \sum_{{\widetilde{d}\left|{m\over \ell}\right.}\atop {(\widetilde{d},\ell)=1}} {{\mu(\widetilde{d})}\over {\widetilde{d}}} 
\buildrel{(1)}\over{=\!=\!=} \1_{\ell|m} {{\mu(\ell)}\over \ell}
                              \sum_{\widetilde{d}\left|{m\over \ell}\right.} {{\mu(\widetilde{d})}\over {\widetilde{d}}} 
\buildrel{(1.14)}\over{=\!=\!=} \1_{\ell|m} {{\mu(\ell)}\over \ell} {{\varphi(m/\ell)}\over {m/\ell}} 
\buildrel{(1)}\over{=\!=\!=} \1_{\ell|m} {{\mu(\ell)\varphi(m)}\over {\varphi(\ell)m}}, 
$$
\par
\noindent
because (in our sums) $(1)$ entails $\mu^2(q)=1$ \& $m|q$ $\Rightarrow \mu^2(m)=1$ whence both : $(\widetilde{d},\ell)=1$ $\forall \widetilde{d}$ dividing $\ell/m$ is implicit AND $\varphi(m/\ell)=\varphi(m)/\varphi(\ell)$, \enspace becomes (since $\ell|m, m|q\,\Rightarrow \ell|q$)
$$
\Win_{\ell}\; \widetilde{{\hbox{\script P}}_{f,g_Q}} 
= {{\mu(\ell)}\over {\varphi(\ell)}}\sum_{{q\le Q}\atop {q\equiv 0(\!\!\bmod \ell)}}\mu(q)\widehat{g_Q}(q)\sum_{{m|q}\atop {m\equiv 0(\!\!\bmod \ell)}}\mu(m)\sum_{p\le N,(p,m)=1}f(p), 
$$
\par
\noindent
which for the moment has an important feature: these two conditions on $q-$sum entail, in particular, that \enspace $\Win_{\ell}\; \widetilde{{\hbox{\script P}}_{f,g_Q}}=0$, $\forall \ell>Q$. Going back to calculations, we split the $p-$sum as
$$
\Win_{\ell}\; \widetilde{{\hbox{\script P}}_{f,g_Q}} 
= {{\mu(\ell)}\over {\varphi(\ell)}}\sum_{{q\le Q}\atop {q\equiv 0(\!\!\bmod \ell)}}\mu(q)\widehat{g_Q}(q)\sum_{{m|q}\atop {m\equiv 0(\!\!\bmod \ell)}}\mu(m)\left(\sum_{p\le N}f(p)-\sum_{p|m}f(p)\right), 
$$
\par
\noindent
so that (hereafter, isolating a term depending on some variables, esp. here $q$ \& $\ell$, these keep their \lq \lq outer sums constraints\rq \rq) 
$$
\sum_{{m|q}\atop {m\equiv 0(\!\!\bmod \ell)}}\mu(m)\buildrel{(1)}\over{=\!=\!=}\mu(\ell)\cdot \1_{\ell|q}\cdot \sum_{\widetilde{m}\left|{q\over {\ell}}\right.}\mu(\widetilde{m})
\buildrel{(1.2)}\over{=\!=\!=} \1_{q=\ell} \enspace \mu(\ell)
\leqno{(4.4)}
$$
\par				
\noindent
implies  
$$
\Win_{\ell}\; \widetilde{{\hbox{\script P}}_{f,g_Q}} 
= {{\widehat{g_Q}(\ell)}\over {\varphi(\ell)}}\sum_{p\le N}f(p)\mu(\ell)
 -{{\mu(\ell)}\over {\varphi(\ell)}}\sum_{{q\le Q}\atop {q\equiv 0(\!\!\bmod \ell)}}\mu(q)\widehat{g_Q}(q)\sum_{p|q}f(p)\sum_{{m|q}\atop {m\equiv 0(\!\!\bmod [p,\ell])}}\mu(m), 
$$
\par
\noindent
which requires (another!) distinction in two cases (we saw in Lemma 4.2), i.e., $p|\ell$ and $p$ doesn't divide $\ell$, (when $p|\ell$, $p\le Q\le N$ is implicit, from $\ell \le Q$ and $\BH$, hereafter)
$$
\Win_{\ell}\; \widetilde{{\hbox{\script P}}_{f,g_Q}} = {{\widehat{g_Q}(\ell)}\over {\varphi(\ell)}}\sum_{p\le N}f(p)\mu(\ell)
 -{{\mu(\ell)}\over {\varphi(\ell)}}\sum_{p|\ell}f(p)\sum_{{q\le Q}\atop {q\equiv 0(\!\!\bmod \ell)}}\mu(q)\widehat{g_Q}(q)\sum_{{m|q}\atop {m\equiv 0(\!\!\bmod \ell)}}\mu(m)
$$
$$
-{{\mu(\ell)}\over {\varphi(\ell)}}\sum_{{q\le Q}\atop {q\equiv 0(\!\!\bmod \ell)}}\mu(q)\widehat{g_Q}(q)\sum_{{p|q}\atop {p\nondiv \ell}}f(p)\sum_{{m|q}\atop {m\equiv 0(\!\!\bmod p\ell)}}\mu(m), 
$$
\par
\noindent
that, thanks to 
$$
\sum_{{m|q}\atop {m\equiv 0(\!\!\bmod p\ell)}}\mu(m)\buildrel{(1)}\over{=\!=\!=}\mu(p\ell)\cdot \1_{(p\ell)|q} \cdot \sum_{\widetilde{m}\left|{q\over {p\ell}}\right.}\mu(\widetilde{m})
\buildrel{(1.2)}\over{=\!=\!=} \1_{q=p\ell} \enspace \mu(p\ell),
$$
\par
\noindent
for all $(p,\ell)=1$, and above $(4.4)$, 
$$
\Win_{\ell}\; \widetilde{{\hbox{\script P}}_{f,g_Q}} = {{\widehat{g_Q}(\ell)}\over {\varphi(\ell)}}\sum_{p\le N}f(p)\mu(\ell)
 -{{\mu(\ell)\widehat{g_Q}(\ell)}\over {\varphi(\ell)}}\sum_{p|\ell}f(p)
  -{{\mu(\ell)}\over {\varphi(\ell)}}\sum_{p\nondiv \ell}f(p)\widehat{g_Q}(p\ell), 
$$
\par
\noindent
where we wish to let appear $c_{\ell}(p)$ instead of $\mu(\ell)$, in the first term, so write: 
$$
\sum_{p\le N}f(p)\mu(\ell)=\sum_{{p\le N}\atop {p\nondiv \ell}}f(p)c_{\ell}(p)+\mu(\ell)\sum_{p|\ell}f(p)
$$
\par
\noindent
whence 
$$
\Win_{\ell}\; \widetilde{{\hbox{\script P}}_{f,g_Q}} = {{\widehat{g_Q}(\ell)}\over {\varphi(\ell)}}\left(\sum_{{p\le N}\atop {p\nondiv \ell}}f(p)c_{\ell}(p)+\mu(\ell)\sum_{p|\ell}f(p)\right)
 -{{\mu(\ell)\widehat{g_Q}(\ell)}\over {\varphi(\ell)}}\sum_{p|\ell}f(p)
  -{{\mu(\ell)}\over {\varphi(\ell)}}\sum_{p\nondiv \ell}f(p)\widehat{g_Q}(p\ell) 
$$
$$
= {{\widehat{g_Q}(\ell)}\over {\varphi(\ell)}}\sum_{{p\le N}\atop {p\nondiv \ell}}f(p)c_{\ell}(p)-{{\mu(\ell)}\over {\varphi(\ell)}}\sum_{p\nondiv \ell}f(p)\widehat{g_Q}(p\ell) 
$$
$$
= {{\widehat{g_Q}(\ell)}\over {\varphi(\ell)}}\sum_{p\le N}f(p)c_{\ell}(p)
 -{{\widehat{g_Q}(\ell)}\over {\varphi(\ell)}}\sum_{p|\ell}f(p)\varphi(\ell){{\mu(\ell/p)}\over {\varphi(\ell/p)}}
   -{{\mu(\ell)}\over {\varphi(\ell)}}\sum_{p\nondiv \ell}f(p)\widehat{g_Q}(p\ell) 
$$
$$
\buildrel{(1)}\over{=\!=\!=} {{\widehat{g_Q}(\ell)}\over {\varphi(\ell)}}\sum_{p\le N}f(p)c_{\ell}(p)
 +\mu(\ell){{\widehat{g_Q}(\ell)}\over {\varphi(\ell)}}\sum_{p|\ell}f(p)\varphi(p)
   -{{\mu(\ell)}\over {\varphi(\ell)}}\sum_{p\nondiv \ell}f(p)\widehat{g_Q}(p\ell), 
$$
\par
\noindent
using : $p|\ell$ $\Rightarrow$ $c_{\ell}(p)=\varphi(\ell)\mu(\ell/p)/\varphi(\ell/p)=-\mu(\ell)\varphi(p)$, thanks to $(1.6)$ and $(1)$.\hfill $\square$

\vfill
\eject

\par				
\smallskip
\par
\noindent {\bf Proof of Th.m 3.4}. From $\widetilde{{\hbox{\script S}}_{f,g_Q}}$ Eratosthenes Transform definition, $(1.4m)$ and $\overline{\chi(-1)}=\chi(-1)$, $\forall d\in \N$, 
$$
\mu^2(d){\hbox{\script S}}_{f,g_Q}'(N,d)=\mu^2(d)\sum_{t|d}\mu\left({d\over t}\right)\sum_{p|t}f(p)p\sum_{q'}\mu(q')\widehat{g_Q}(pq')\sum_{m'|q'}{{\mu(m')m'}\over {\varphi(m')}}\sum_{\chi\atop {(\!\!\bmod m')}}\chi(-1)\chi(t/p)
$$
$$
=\sum_{p|d}f(p)p\sum_{q'}\mu(q')\widehat{g_Q}(pq')\sum_{m'|q'}{{\mu(m')m'}\over {\varphi(m')}}\sum_{\chi\atop {(\!\!\bmod m')}}\chi(-1)\mu(d)\sum_{{t|d}\atop {t\equiv 0(\!\!\bmod p)}}\mu(t)\chi(t/p) = (\hbox{\rm set}\enspace t':=t/p) =
$$
$$
=\mu(d)\sum_{p|d,p\le Q}f(p)p\sum_{q'}\mu(q')\widehat{g_Q}(pq')\sum_{m'|q'}{{\mu(m')m'}\over {\varphi(m')}}\sum_{\chi\atop {(\!\!\bmod m')}}\chi(-1)\left(-\sum_{t'\left|{d\over p}\right.}\mu(t')\chi(t')\right), 
$$
\par
\noindent
whence,by $\S2.1$ formula and Wintner $\ell-$th coefficient definition,  
$$
\Win_{\ell} \enspace \widetilde{{\hbox{\script S}}_{f,g_Q}}=
-\sum_{p\le Q}f(p)p\sum_{q'}\mu(q')\widehat{g_Q}(pq')\sum_{m'|q'}{{\mu(m')m'}\over {\varphi(m')}}\sum_{\chi\atop {(\!\!\bmod m')}}\chi(-1)\Sigma_p(\chi,\ell,m'), 
$$
\par
\noindent
from Lemma 4.2 definition, here also for $\chi=\chi_0(\bmod \;m')$, whence Pinch Lemma 4.1 for $\chi_0(\bmod \;m')$, here, gives $(d/p)|m'$, in the $d-$sum: 
$$
\Win_{\ell} \enspace \widetilde{{\hbox{\script S}}_{f,g_Q}}=
-\sum_{p\le Q}f(p)p\sum_{q'}\mu(q')\widehat{g_Q}(pq')\sum_{m'|q'}{{\mu(m')m'}\over {\varphi(m')}}\sum_{{d\equiv 0(\!\!\bmod \ell)}\atop {{d\equiv 0(\!\!\bmod p)}\atop {(d/p)|m'}}}{{\mu(d)}\over d}, 
$$
\par
\noindent
where (recall, isolating a term, its variables, esp. here $p,\ell,m'$, keep their \lq \lq outer sums constraints\rq \rq) 
$$
\Sigma_p(\chi_0,\ell,m')=\sum_{{d\equiv 0(\!\!\bmod \ell)}\atop {{d\equiv 0(\!\!\bmod p)}\atop {(d/p)|m'}}}{{\mu(d)}\over d}
= (\hbox{\rm set} \enspace d':=d/p \enspace \hbox{\rm here}) 
= \sum_{{d'|m'}\atop {pd'\equiv 0(\!\!\bmod \ell)}}{{\mu(pd')}\over {pd'}}
= -{1\over p}\sum_{{(d',p)=1}\atop {{d'|m'}\atop {pd'\equiv 0(\!\!\bmod \ell)}}}{{\mu(d')}\over {d'}}
$$
\par
\noindent
needing again the splitting $p|\ell$ aut $p$ doesn't divide $\ell$ : 
$$
\Sigma_p(\chi_0,\ell,m')=-{1\over p}\left(\1_{p|\ell}
\sum_{
{(d',p)=1}\atop {{d'|m'}\atop {d'\equiv 0\left(\!\!\bmod {{\ell}\over p}\right)}}
}
{{\mu(d')}\over {d'}}
+\1_{p\nondiv \ell}\sum_{{(d',p)=1}\atop {{d'|m'}\atop {d'\equiv 0(\!\!\bmod \ell)}}}{{\mu(d')}\over {d'}}\right)
$$
\par
\noindent
and, resp. \enspace $d'={{\ell}\over p}\widetilde{d}$ \enspace for $p|\ell$, resp., \enspace $d'=\ell \widetilde{d}$ \enspace otherwise, give 
$$
\Sigma_p(\chi_0,\ell,m')=-{1\over p}\left(
\1_{p|\ell}\cdot \1_{{{\ell}\over p}\left|m'\right.}\cdot {{\mu(\ell/p)}\over {\varphi(\ell/p)}}
 \cdot \sum_{{(\widetilde{d},\ell)=1}\atop {\widetilde{d}\left|{{m'}\over {\ell/p}}\right.}}{{\mu(\widetilde{d})}\over {\widetilde{d}}}
+\1_{p\nondiv \ell}\cdot \1_{\ell|m'}\cdot {{\mu(\ell)}\over {\ell}}\sum_{{(\widetilde{d},p\ell)=1}\atop {\widetilde{d}\left|{{m'}\over {\ell}}\right.}}{{\mu(\widetilde{d})}\over {\widetilde{d}}}\right). 
$$
\par
\noindent
This needs both 
$$
\sum_{{(\widetilde{d},\ell)=1}\atop {\widetilde{d}\left|{{m'}\over {\ell/p}}\right.}}
{{\mu(\widetilde{d})}\over {\widetilde{d}}}
\buildrel{(1.3)}\over{=\!=\!=} 
\prod_{{\widetilde{p}\left|{{m'}\over {\ell/p}}\right.}\atop {\widetilde{p}\nondiv \ell}}\left(1-{1\over {\widetilde{p}}}\right)
\buildrel{(A)}\over{=\!=\!=} \prod_{\widetilde{p}\left|{{m'}\over {\ell/p}}\right.}\left(1-{1\over {\widetilde{p}}}\right)
\buildrel{\mu(m')\neq 0}\over{=\!=\!=\!=\!=} {{\varphi(m')}\over {m'}}\cdot {{\ell/p}\over {\varphi(\ell/p)}}
$$
\par				
\noindent
and 
$$
\sum_{{(\widetilde{d},p\ell)=1}\atop {\widetilde{d}\left|{{m'}\over {\ell}}\right.}}
{{\mu(\widetilde{d})}\over {\widetilde{d}}}
\buildrel{(1.3)}\over{=\!=\!=} 
\prod_{{\widetilde{p}\left|{{m'}\over {\ell}}\right.}\atop {\widetilde{p}\nondiv (p\ell)}}\left(1-{1\over {\widetilde{p}}}\right)
\buildrel{(B)}\over{=\!=\!=} \prod_{\widetilde{p}\left|{{m'}\over {\ell}}\right.}\left(1-{1\over {\widetilde{p}}}\right)
\buildrel{\mu(m')\neq 0}\over{=\!=\!=\!=\!=} {{\varphi(m')}\over {m'}}\cdot {{\ell}\over {\varphi(\ell)}}, 
$$
\par
\noindent
where we justify $(A)$ from : IF $\exists P\in \Primes$ with $P|\ell$ \& $P|{{m'}\over {\ell/p}}$ $\Rightarrow$ $P^2|(m'p)$, $m'|q'$, THEN $P^2|(pq')$ \enspace $\assurdo$ (abbreviating: absurd), from Hypothesis $(1)$; we may justify $(B)$ similarly, again from Hypothesis $(1)$: IF $\exists P\in \Primes$ with $P|(p\ell)$ \& $P|{{m'}\over {\ell}}$ $\Rightarrow$ $P^2|(m'p)$, THEN as above $\assurdo$. 
\par
\noindent
Thus
$$
\Sigma_p(\chi_0,\ell,m')=-{1\over p}\left(
\1_{p|\ell}\cdot \1_{{{\ell}\over p}\left|m'\right.}\cdot {{\mu(\ell/p)}\over {\varphi(\ell/p)}}
+\1_{p\nondiv \ell}\cdot \1_{\ell|m'}\cdot {{\mu(\ell)}\over {\varphi(\ell)}}\right){{\varphi(m')}\over {m'}}
$$
$$
=-{1\over p}\cdot {{\varphi(m')}\over {m'}}\left(
\1_{p\nondiv \ell}\cdot \1_{\ell|m'}
-\1_{p|\ell}\cdot \1_{{{\ell}\over p}\left|m'\right.}\cdot \varphi(p)
\right){{\mu(\ell)}\over {\varphi(\ell)}},
\enspace
\forall \ell \enspace \hbox{\rm square-free}. 
$$
\par
\noindent
Hence,
$$
\Win_{\ell} \enspace \widetilde{{\hbox{\script S}}_{f,g_Q}}=
{{\mu(\ell)}\over {\varphi(\ell)}}\sum_{p\le Q}f(p)\sum_{q'}\mu(q')\widehat{g_Q}(pq')\sum_{m'|q'}\mu(m')
\left(
\1_{p\nondiv \ell}\cdot \1_{\ell|m'}
-\1_{p|\ell}\cdot \1_{{{\ell}\over p}\left|m'\right.}\cdot \varphi(p)
\right)  
$$
$$
=
{{\mu(\ell)}\over {\varphi(\ell)}}
\left(
\sum_{{p\nondiv \ell}\atop {p\le Q}}f(p)\sum_{q'\equiv 0(\!\!\bmod \ell)}\mu(q')\widehat{g_Q}(pq')\sum_{{m'|q'}\atop {m'\equiv 0(\!\!\bmod \ell)}}\mu(m')
\right.
$$
$$
\left.
-\sum_{p|\ell}f(p)\varphi(p)\sum_{q'\equiv 0\left(\!\!\bmod {{\ell}\over p}\right)}\mu(q')\widehat{g_Q}(pq')\sum_{{m'|q'}\atop {m'\equiv 0\left(\!\!\bmod {{\ell}\over p}\right)}}\mu(m')
\right). 
$$
\par
\noindent
We need only 
$$
\sum_{{m'|q'}\atop {m'\equiv 0(\!\!\bmod \ell)}}\mu(m')
\buildrel{(C)}\over{=\!=\!=} \mu(\ell)\sum_{\widetilde{m}\left|{{q'}\over {\ell}}\right.}\mu(\widetilde{m})
\buildrel{(1.2)}\over{=\!=\!=} \mu(\ell)\cdot \1_{q'=\ell}
\quad\hbox{\stampatello and}
$$
$$
\sum_{{m'|q'}\atop {m'\equiv 0\left(\!\!\bmod {{\ell}\over p}\right)}}\mu(m')
\buildrel{(D)}\over{=\!=\!=} \mu\left({{\ell}\over p}\right)\sum_{\widetilde{m}\left|{{q'}\over {\ell/p}}\right.}\mu(\widetilde{m})
\buildrel{(1.2)}\over{=\!=\!=} \mu\left({{\ell}\over p}\right)\cdot \1_{q'={{\ell}\over p}}, 
$$
\par
\noindent
where in $(C)$: {\stampatello condition $(\widetilde{m},\ell)=1$ is for free}: IF $\exists P|\ell$ \& $P|(q'/\ell)$ $\Rightarrow$ $P^2|q'$ THEN from Hp. $(1)$, $\assurdo$
\par
\noindent
and in $(D)$: {\stampatello condition $(\widetilde{m},\ell/p)=1$ holds o.w.}: IF $\exists P|(\ell/p)$ \& $P{\displaystyle \left|{{q'}\over {\ell/p}}\right. }$ $\Rightarrow$ $P^2|q'$ THEN as above $\assurdo$. 
\par
\noindent
Thus
$$
\Win_{\ell} \enspace \widetilde{{\hbox{\script S}}_{f,g_Q}}=
{{\mu(\ell)}\over {\varphi(\ell)}}
\left(\sum_{{p\nondiv \ell}\atop {p\le Q}}f(p)\mu^2(\ell)\widehat{g_Q}(p\ell)
 -\sum_{p|\ell}f(p)\varphi(p)\mu^2\left({{\ell}\over p}\right)\widehat{g_Q}(\ell)
\right) 
$$
$$
\buildrel{(1)}\over{=\!=\!=} -{{\mu(\ell)}\over {\varphi(\ell)}}\widehat{g_Q}(\ell)\sum_{p|\ell}f(p)\varphi(p) + {{\mu(\ell)}\over {\varphi(\ell)}}\sum_{{p\nondiv \ell}\atop {p\le Q}}f(p)\widehat{g_Q}(p\ell).
$$
\hfill $\square$

\vfill
\eject

\par				
\noindent{\bf 4.5. Three Hypotheses Correlations: an assumption more gives H-L Conjecture for $\sf a=2k$}
\bigskip
\par
\noindent
We begin, towards: proving conditionally H-L Conjecture for \sf shifts $a=2k$, from calculating the $P-$smooth Wintner coefficients of Three-Hypotheses Correlations from an adaption of Lemma 4.2 to these, say, $P-$local coefficients; they are defined, in general, for a fixed $F:\N \rightarrow \C$, as, once FIXED $P\in \Primes$ (recall the $P-$smooth naturals, $(P)$, from page$5$) : 
$$
\forall \ell \in (P),\quad \WinPellF \defineq \sum_{{d\in (P)}\atop {d\equiv 0(\!\!\bmod \ell)}}{{F'(d)}\over d}. 
$$
\par
\noindent
Of course these need not to exist, as the series might not converge (and we require convergence, for all $\ell \in (P)$ involved); trivially, of course, it converges, to \enspace $0$, for all \enspace $\ell \in \N \backslash (P)$ \hfill as the series is an empty sum. 
\smallskip
\par
\noindent
Recall (see page 8) a superscript \lq \lq $\flat$\rq \rq\enspace restricts to \sf, both for sums: $\mathop{{\sum}^{\flat}}$, and for elements in $S^{\flat}$, $\forall S\subseteq \N$; in particular, $S=(P)$, the $P-$smooth set, gives rise to $\Pb \defineq \{n\in \N : n \enspace \sf \hbox{\stampatello and} \enspace n\in (P)\}$. \hfill Same for $\N^{\flat}$. 
\smallskip
\par
Instead, fixed any $P\in \Primes$, for $\widetilde{F}$ the $\ell-$th $P-$local Wintner coefficient exists (recall characterization of Eratosthenes Transform, see $\S1$, of $\IPP$ functions in $\S2.1$), $\forall \ell \in \N$, vanishing outside $\Pb$, this time: 
$$
\forall \ell \in \Pb,\quad \WinPellFtilde \defineq \sum_{{d\in (P)}\atop {d\equiv 0(\!\!\bmod \ell)}}{{\mu^2(d)F'(d)}\over d}
 \EqByDef \sum_{{d\in \Pb}\atop {d\equiv 0(\!\!\bmod \ell)}}{{F'(d)}\over d}. 
$$
\par
\noindent
We have : two advantages, uniformly whatever $F$ we choose, studying $\widetilde{F}$ in place of our $F$; the first, these coefficients above for $\widetilde{F}$ exist (again, uniformly w.r.t. $F:\N \rightarrow \C$ we choose), since $\Pb$, once fixed $P\in \Primes$, is finite (as $|\Pb|=2^{\pi(P)}$, $\forall P\in \Primes$), so finite support, uniformly w.r.t. $\ell \in \N$, of above sums gives convergence to ALL $P-$local coefficients for $\widetilde{F}$; second, $|\Pb|<\infty$ entails a {\it finite Ramanujan expansion} for $F(a)$, once fixed $a\in \Pb$. Compare [C2] for this, in particular; and for  the properties of these $P-$local coefficients, both for $F$ and for $\widetilde{F}$. 
\par 
Here, we'll confine to $F$ to be a three-hypotheses Correlation, for which particular $F$ we get: above $P-$Wintner coefficients for $\widetilde{F}$, in the limit over $P$ running in primes, converge to the Wintner coefficients for $\widetilde{F}$. This is the content of Theorem 3.5 above: for the Proof we need following Lemma 4.3, coming from Lemma 4.2 arguments. 
\medskip
\par
We calculate the $P-$limit for a kind of $P-$local versions of $\Sigma(\chi,\ell,m)$ \& $\Sigma_p(\chi,\ell,m')$ (in Lemma 4.2). 
\smallskip
\par
\noindent {\bf Lemma 4.3}. ({\stampatello $P$-local divisors version of Lemma 4.2})
\par
\noindent
{\it Let } $p\in \Primes$, $P\in \Primes$, $\ell\in \N$. {\it Fix } $m\in \N, m>2$ {\it and fix any non-principal Dirichlet character } $\chi\neq \chi_0(\bmod \;m)$. {\it Define}
$$
\Sigma^{(P)}(\chi,\ell,m)\defineq \sum_{{d\in (P)}\atop{d\equiv 0(\!\!\bmod \ell)}}{{\mu(d)}\over d}\sum_{t|d}\mu(t)\chi(t)
 \EqByDef \Win_{\ell}^{(P)}\; \widetilde{\chi} , 
$$
\par
\noindent
{\it while, fixed } $m'\in \N,m'>2$ {\it and fixed any non-principal Dirichlet character } $\chi\neq \chi_0(\bmod \;m')$, {\it define}
$$
\Sigma^{(P)}_p(\chi,\ell,m')\defineq 
\sum_{
{
{d\in (P)}
\atop
{{d\equiv 0(\!\!\bmod \ell)}
\atop {d\equiv 0(\!\!\bmod p)}}}
}
{{\mu(d)}\over d}\sum_{t'\left|{d\over p}\right.}\mu(t')\chi(t'). 
$$
\par
\noindent
{\it Then}
$$
\lim_P\enspace \Sigma^{(P)}(\chi,\ell,m)=\lim_P\enspace \Sigma^{(P)}_p(\chi,\ell,m')=0, 
\quad \hbox {\stampatello uniformly}\enspace \forall \ell \in \N. 
$$

\vfill

\par
\noindent {\bf Remark 3.} The Proof for \enspace $\Sigma^{(P)}_p(\chi,\ell,m')$ \enspace follows the {\stampatello same lines} for \enspace $\Sigma_p(\chi,\ell,m')$, $\underline{\hbox{\stampatello after}}$ {\stampatello distinguishing two cases}: $p|\ell$ {\stampatello aut} $p\not | \ell$, mimicking following argument for $\Sigma^{(P)}(\chi,\ell,m)$.\hfill $\diamond$  

\eject

\par				
\noindent{\bf Proof}. We wish to prove, since $\Sigma^{(P)}(\chi,\ell,m)\EqByDef \Win_{\ell}^{(P)}\; \widetilde{\chi}$ , 
$$
\chi \neq \chi_0(\!\bmod\; m)
\quad \Rightarrow \quad 
\Win_{\ell}^{(P)}\; \widetilde{\chi}\ll_m {1\over {\log P}}, 
\enspace \hbox{\stampatello uniformly} \enspace \forall \ell \in \N,
\leqno{(\ast)_{\ell}}
$$
\par
\noindent
whence : \enspace $\lim_P \Win_{\ell}^{(P)}\; \widetilde{\chi}=0$, {\stampatello unif.}$^{\hbox{\stampatello ly}}$ $\forall \ell \in \N$. We may assume $\ell \in \Pb$, since otherwise LHS vanishes. Hereafter, we have to fix \enspace $\ell \in \Pb$, so
$$
\Win_{\ell}^{(P)}\; \widetilde{\chi}\EqByDef \sum_{{d\in \Pb}\atop {d\equiv 0(\!\!\bmod \ell)}}{{\chi'(d)}\over d}
 ={{\chi'(\ell)}\over {\ell}} \sum_{{n\in \Pb}\atop {(n,\ell)=1}}{{\chi'(n)}\over n}
  ={{\chi'(\ell)}\over {\ell}} \prod_{{p\le P}\atop {p\nondiv \ell}}\left(1+{{\chi'(p)}\over p}\right), 
$$
\par
\noindent
using $(1.3)$ above; by definition, $\chi'(p)=\chi(p)-1$, $\forall p\in \Primes$ and for any Dirichlet character $\chi$; ignoring the prime $p=2$ (compare Lemma 4.2 Proof), that for $\ell$ odd makes LHS vanish, when the (necessarily) non-principal $\chi(\bmod \; m)$ has $\chi(2)=-1$, we have to show, say, in order to prove $(\ast)_\ell$ $\!$: 
$$
\chi \neq \chi_0(\!\bmod\; m)
\thinspace \Rightarrow \thinspace 
W:={{\chi'(\ell)}\over {\ell}} \prod_{2<p\le P}\left(1-{1\over p}+{{\chi(p)}\over p}\right) \prod_{{p|\ell}\atop {p>2}}\left(1+{{\chi'(p)}\over p}\right)^{-1}
\ll_m {1\over {\log P}}, 
\enspace \hbox{\stampatello unif.}^{\hbox{\stampatello ly}} \enspace \forall \ell \in \N;
$$
\par
\noindent
we start, then, from estimating the $\ell-$dependent two factors: the former is
$$
{{\chi'(\ell)}\over {\ell}}\EqByDef {1\over {\ell}}\sum_{d|\ell}\chi(d)\mu(\ell/d)
\ll {{d(\ell)}\over {\ell}}, 
$$
\par
\noindent
where the number of $\ell$ divisors is $d(\ell)=2^{\omega(\ell)}$, being $\ell$ square-free; while the second's
$$
\prod_{{p|\ell}\atop {p>2}}\left(1+{{\chi'(p)}\over p}\right)^{-1}=\prod_{{p|\ell}\atop {p>2}}\left({{p+\chi(p)-1}\over p}\right)^{-1}
 =\prod_{{p|\ell}\atop {p>2}}\left({p\over {p-1+\chi(p)}}\right)
  \ll \prod_{{p|\ell}\atop {p>2}}{p\over {p-2}} ,
$$
\par
\noindent
uniformly for all Dirichlet characters (having, for all primes $p$, $|p-1+\chi(p)|\ge p-1-|\chi(p)|\ge p-2$ : recall $p>2$ in above product), making
$$
\prod_{{p|\ell}\atop {p>2}}\left(1+{{\chi'(p)}\over p}\right)^{-1}\ll \prod_{{p|\ell}\atop {p>2}}\left(1+{2\over {p-2}}\right)
 \ll \prod_{{p|\ell}\atop {p\ge 5}}\left(1+{2\over 3}\right)
  \ll \prod_{p|\ell}2 = 2^{\omega(\ell)}; 
$$
\par
\noindent
these two factors' bounds give, using $d(\ell)\cdot 2^{\omega(\ell)}=4^{\omega(\ell)}$, for \sf $\ell$, recalling $W$ local definition above: 
$$
W\ll \exp\left({\sum_{2<p\le P}\log\left(1-{1\over p}+{{\chi(p)}\over p}\right)}\right) 
 \ll \exp\left(\sum_{2<p\le P}{{\chi(p)-1}\over p}\right),
$$
\par
\noindent
since : 
$$
{{4^{\omega(\ell)}\over {\ell}}}=\prod_{p|\ell}{4\over p}\le {{4^2}\over 6}\prod_{{p|\ell}\atop {p\ge 5}}{4\over p}\ll 1, 
$$
\par
\noindent
with absolute constant $8/3$. Using $(4.3)$, as we did in Lemma 4.2 Proof (this time $x=P\to \infty$), and 
$$
\sum_{p\le x}{1\over p}=\log\log x+O(1), 
$$
\par
\noindent
as $x\to \infty$ [T], with $x=P$ again, above bounds for $W$ are entailing $(\ast)_{\ell}$ above.\hfill $\square$ 

\vfill
\eject

\par				
\noindent
Following {\stampatello Lemma 4.3 in place of above Lemma 4.2}, we prove very easily {\stampatello Theorem 3.5}.
\smallskip
\par
\noindent{\bf Proof of Theorem 3.5}. Like in above Proof of {\stampatello Theorem 3.3}, we follow {\stampatello Theorem 3.2} to express Correlation's {\stampatello primary part} in terms of $\chi(\!\bmod\; m)$ : while {\stampatello Th.m 3.3} has {\stampatello vanishing contribution} from the $\chi\neq \chi_0(\!\bmod\; m)$, {\stampatello Th.m 3.5} has {\stampatello infinitesimal contribution} but this is true {\stampatello uniformly} $\forall \ell\in \N$, as $P\to \infty$, by Lemma 4.3. {\stampatello Following th.m 3.3} Proof, we arrive at
$$
\lim_P\, \Win_{\ell}^{(P)}\; \widetilde{{\hbox{\script P}}_{f,g_Q}} = \sum_{q\le Q}\mu(q)\widehat{g_Q}(q)\sum_{m|q}{{\mu(m)m}\over {\varphi(m)}}\sum_{p\le N,(p,m)=1}f(p)\,\hbox{\script L}\, ,
\leqno{(\ast)_{\hbox{\script L}}}
$$
\par
\noindent
say, where \enspace
$$
\hbox{\script L}:=\lim_P\, \sum_{{{d\in (P)}\atop {{d\equiv 0(\!\!\bmod \ell)}\atop {d|m}}}}{{\mu(d)}\over d}
 =\sum_{{{d\equiv 0(\!\!\bmod \ell)}\atop {d|m}}}{{\mu(d)}\over d}, 
\enspace \hbox{\stampatello uniformly} \enspace
\forall \ell \in \N, 
$$ 
\par
\noindent
from : $d|m\, \Rightarrow \, d|q\, \Rightarrow \, d\le Q$, using $m,q$ constraints in $(\ast)_{\hbox{\script L}}$ above. Hence,
\smallskip
\par
\leftline{
$
\hfill {\displaystyle \lim_P\; \Win^{(P)}_{\ell} \widetilde{{\hbox{\script P}}_{f,g_Q}} = \Win_{\ell} \; \widetilde{{\hbox{\script P}}_{f,g_Q}} ,
\enspace \hbox{\stampatello uniformly} \enspace
\forall \ell \in \N. }\hfill {\rm QED}
$
}
\par
\noindent
{\stampatello Summarizing}, this holds: {\stampatello from} the {\stampatello non-principal} $\chi\neq \chi_0(\!\bmod\; m)$ {\stampatello infinitesimal contribution} (as \enspace $P\to \infty$, $\hbox{\stampatello unif.}^{\hbox{\stampatello ly}}$ \enspace $\forall \ell \in \N$), {\stampatello while, in the limit} $P\to \infty$, the {\stampatello principal} $\chi=\chi_0(\!\bmod\; m)$ {\stampatello contributes} to $P-$local Wintner $\ell-$th coefficient {\stampatello with the Wintner $\ell-$th coefficient.} 
\smallskip
\par
\noindent
{\stampatello The same argument} (with small changes) {\stampatello works for} $\widetilde{{\hbox{\script S}}_{f,g_Q}}$ : we {\stampatello follow th.m 3.4} Proof \& , after Lemma 4.3 application, we take \thinspace $\chi=\chi_0(\!\bmod\; m')$ \thinspace and we {\stampatello adapt \thinspace $(\ast)_{\hbox{\script L}}$ \thinspace above}, using: $\left.{d\over p}\right| m'\thinspace \Rightarrow\thinspace \left.{d\over p}\right| q'\thinspace \Rightarrow\thinspace d\le Q$.\hfill $\square$  

\bigskip

\par
We go into the \lq \lq deep\rq \rq, so-to-speak: the link in between the $P-$Wintner coefficients \& the Ramanujan coefficients, both for a very general $\IPP$ function, namely $\widetilde{F}$, for general $F:\N \rightarrow \C$. 
\medskip
\par
The technical definitions we need in next Proof will be given inside, as they are needed, at right time. 
\smallskip
\par
\noindent {\bf Proof of Theorem 3.6}. We begin by an important distinction: next implication is NOT an equivalence
$$
\sum_{q\;\cf}|G(q)|<\infty 
\quad \Rightarrow \quad 
\forall a\in \N^{\flat},\enspace \sum_{q\;\cf}\left|G(q)c_q(a)\right|=\sum_{q=1}^{\infty}\left|G(q)c_q(a)\right|<\infty 
$$
\par
\noindent
and it holds, as the absolute convergence of Ramanujan Coefficients (for $\widetilde{F}$) on cube-free naturals implies the absolute convergence of Ramanujan Expansion (for $\widetilde{F}$), on all $a\in \N^{\flat}$. This follows from two properties, both regarding the Ramanujan sums, the first is at page 15: {\stampatello Ramanujan vertical limit}, abbrev. $\Rvl$, which implies, for \sf arguments $a\in \N$, that Ramanujan Expansion has a cube-free supported Ramanujan Coefficient; by $\Rvl$, we may {\bf assume henceforth}: {\stampatello $G$ is \cf supported}. The second property comes from Lemma A.1 in [CM]: $|c_q(a)|\le (a,q)$, $\forall q\in \N$, $\forall a\in \Z$ , whence trivially $|c_q(a)|\le a$, $\forall a,q\in \N$. Together, these two properties prove above implication. Converse implication doesn't hold, as we don't have a general lower bound for $|c_q(a)|$, at least, not for \cf $q\in \N$ (while $q$ \sf $\Rightarrow$ $|c_q(a)|\ge 1$, $\forall a\in \N$ [C2]). Hence, the property we require to Ramanujan Expansion for our $\widetilde{F}$, on \sf $a\in \N$, is NOT ONLY the ABSOLUTE CONVERGENCE OF R.E. but of its Ramanujan Coefficients (that above implication proves to be a stronger assumption and this is the reason why we wrote both hypotheses in the statement, inside \lq \lq following property\rq \rq). 
\par
We need to show $(3.2)$ and for this we fix $P\in \Primes$ henceforth and we {\stampatello start from}:
$$
F_{\Pb}(a)\defineq \sum_{{d\in \Pb}\atop {d|a}}F'(d), 
\enspace 
\forall a\in \N
$$
\par
\noindent
which is, for technical reasons, needed here. (Technical means: M\"obius inversion, that starts from a function defined on all $a\in \N$.) See that, defining the $\Pb-$component of a natural argument $a\in \N$, written $a_{\Pb}$ : 
$$
a_{\Pb}\defineq \prod_{{p\in (P)}\atop {p|a}}p
 \EqByDef \prod_{{p\le P}\atop {p|a}}p
\quad \Rightarrow \quad 
F_{\Pb}(a)\EqByDef F\left(a_{\Pb}\right), 
\enspace 
\forall a\in \N. 
$$
\par
\noindent
Furthermore, we have \enspace $F_{\Pb}(a)\EqByDef F\left(a_{\Pb}\right)\EqByDef \widetilde{F}\left(a_{\Pb}\right)$, $\forall a\in \N$, which we use at once. 

\vfill
\eject

\par				
\noindent
{\stampatello Thus} (from $\widetilde{F}$ property in the statement)
$$
F_{\Pb}(a)=\sum_{q=1}^{\infty}G(q)c_q(a_{\Pb}),
\enspace \forall a\in \N,
$$
\par
\noindent
by above $F_{\Pb}$ definition; write \thinspace $q=rs$, \thinspace with \thinspace $s:=q_{(P)}\defineq \prod_{p\le P}p^{v_p(q)}$ \thinspace and  \thinspace $r:=q_{)P(}\defineq \prod_{p>P}p^{v_p(q)}$ \hfill (say, 
\smallskip
\par
\noindent
the $s${\stampatello mooth} \& the $r${\stampatello ough} components of $q$, with $s\in (P)$, in the $P-${\stampatello smooth} and $r\in )P($, in the $P-${\stampatello sifted}). 
\par
\noindent
This and: $c_q(a)$ {\stampatello multiplicative} w.r.t. $q\in \N$ entail $c_q\left(a_{\Pb}\right)=\mu(r)c_s\left(a_{\Pb}\right)$, from $c_r\left(a_{\Pb}\right)=\mu(r)$, because, in $(1.6)$, $(r,a_{\Pb})=1$ {\stampatello by definitions} above; whence
$$
F_{\Pb}(a)=\sum_{s\in (P)}c_s\left(a_{\Pb}\right)\sum_{r\in )P(}\mu(r)G(sr), 
\enspace \forall a\in \N. 
$$
\par
\noindent
{\stampatello Hereafter} all sums exchanges are justified, thanks to the Absolute Convergence of above Ramanujan Expansion discussed earlier.  Furthermore, here, $sr$ is \cf as $r\in )P($ is \sf (thanks to $\mu(r)-${\stampatello twist}) and, applying $\Rvl$ above, we have \cf $s\in (P)$; whence 
$$
F_{\Pb}(a)=\sum_{s\in (P)}c_s\left(a_{\Pb}\right)\sum_{r\in )P(}\mu(r)\Gcf(sr), 
\enspace \forall a\in \N, 
$$
\par
\noindent
hereafter writing $\Gcf$ to enhance: $G$ is \cf supported. Applying Kluyver's $(1.7)$, then exchanging sums, 
$$
F_{\Pb}(a)=\sum_{{d|a}\atop {d\in \Pb}}d\sum_{K\in (P)}\mu(K)\sum_{r\in )P(}\mu(r)\Gcf(dKr), 
\enspace \forall a\in \N,
$$
\par
\noindent
whence (finally!) applying M\"obius inversion soon after $(1.2)$ (since above formula is $\forall a\in \N$), together with above $F_{\Pb}$ definition, 
$$
\forall d\in \Pb, \enspace F'(d)=d\sum_{K\in (P)}\mu(K)\sum_{r\in )P(}\mu(r)\Gcf(dKr). 
$$
\par
\noindent
{\stampatello Thus}
$$
\Win_{\ell}^{(P)}\; \widetilde{F} \EqByDef \sum_{{d\in \Pb}\atop {d\equiv 0(\!\!\bmod\, \ell)}}\sum_{K\in (P)}\mu(K)\sum_{r\in )P(}\mu(r)\Gcf(drK),
\quad \forall \ell \in \Pb. 
$$
\par
\noindent
Here $d$ \sf, writing $d=\ell m$, {\stampatello entails a factor} $\mu^2(\ell m)$: 
$$
\forall \ell \in \Pb, 
\quad
\Win_{\ell}^{(P)}\; \widetilde{F}=\sum_{m\in \Pb}\mu^2(\ell m)\sum_{K\in (P)}\mu(K)\sum_{r\in )P(}\mu(r)\Gcf(\ell rKm), 
\enspace \hbox{\stampatello namely}:
$$
$$
\forall \ell \in \Pb, 
\enspace 
\Win_{\ell}^{(P)}\; \widetilde{F}=\sum_{{m\in \Pb}\atop {(m,\ell)=1}}\sum_{K\in (P)}\mu(K)\sum_{r\in )P(}\mu(r)\Gcf(\ell rKm)
 =\sum_{n\in (P)}\sum_{r\in )P(}\mu(r)\Gcf(\ell rn)S_{\ell}(n), 
$$
\par
\noindent
{\stampatello say, defining}:
$$
\forall \ell \in \Pb,\enspace \forall n\in (P),
\quad
S_{\ell}(n)\defineq \sum_{{m|n}\atop {(m,\ell)=1}}\mu^2(m)\mu\left({n\over m}\right). 
$$
\par
\noindent
Since \enspace $\ell rn$ \enspace has to be \cf \& \cf naturals are divisor-closed, we have the {\stampatello constraint } $n$ \cf : use $(2.1)$, with $n=st^2$ in place of $q=uv^2$ \& $u=s$, $v=t$ both \sf and $(s,t)=1$, to get, with $*-$sums meaning summation variables coprime to $\ell$, 
$$
S_{\ell}\left(st^2\right)=\starsum_{m_s|s}\starsum_{m_t|t}\mu\left({{st^2}\over {m_s m_t}}\right)
 =\starsum_{m_s|s}\starsum_{m_t|t}\mu\left({s\over {m_s}}\cdot {t\over {m_t}}\cdot t\right), 
\qquad
\hbox{\stampatello using\enspace twice}\enspace (s,t)=1. 
$$
\par
\noindent
Within the M\"obius function above, $t$ and $t/m_t$ have to be coprime (o.w., their product becomes not \sf and this makes $\mu$ vanish): hence, we have the only possibility\enspace $m_t=t$; using\enspace $(s,t)=1$\enspace again,
$$
S_{\ell}\left(st^2\right)=\1_{(t,\ell)=1}\mu(t)\starsum_{m_s|s}\mu\left({s\over {m_s}}\right)
 =\1_{(t,\ell)=1}\mu(t)\1_{s|\ell}\mu(s), 
\enspace 
\hbox{\stampatello by\thinspace next\thinspace Lemma}\thinspace 5.1\thinspace \&\thinspace \hbox{\stampatello Remark}\thinspace 4. 
$$

\vfill
\eject

\par				
\noindent
Back to our $P-$local $\ell-$th Wintner coefficient for $\widetilde{F}$, 
$$
\forall \ell \in \Pb, 
\enspace 
\Win_{\ell}^{(P)}\; \widetilde{F}=\sum_{{s\in (P)}\atop {s|\ell}}\mu(s)\sum_{{t\in (P)}\atop {(t,\ell)=1}}\mu(t)\sum_{r\in )P(}\mu(r)\Gcf(\ell rst^2),
$$
\par
\noindent
where: $s$ \& $t$ are \sf by $\mu(s)$ \& $\mu(t)$\enspace {\stampatello twisting}, $(s,t)=1$ {\stampatello from } $s|\ell$ \& $(t,\ell)=1$. Also, since $(P)$ is a divisor-closed set, $s|\ell$ \& $\ell \in \Pb$ $\Rightarrow $ $s\in (P)$: \enspace $\forall \ell \in \Pb$, 
$$
\Win_{\ell}^{(P)}\; \widetilde{F}=\sum_{s|\ell}\mu(s)\sum_{{t\in (P)}\atop {(t,\ell)=1}}\mu(t)\sum_{r\in )P(}\mu(r)G(\ell rst^2),
$$
\par
\noindent
{\stampatello i.e., $(3.2)$ above}.\hfill QED
\par
\noindent
We start, now, proving next formula, for $\lim_P$, with a {\stampatello fundamental remark} : in just proved $(3.2)$ we have the $G-$argument, $\ell rst^2$, still \cf as $\ell \in \Pb$ is fixed and summation variables $r,s,t$ constraints, of course, entail this. The only thing changing is the {\it unique representation} given in $(2.1)$ applied here to: 
$$
\ell rst^2 = \left({{\ell}\over s}\cdot r\right)\cdot \left(st\right)^2
 \buildrel{\hbox{\stampatello say}}\over{=\!=\!=} uv^2, 
\enspace u\defineq {{\ell}\over s}\cdot r,\enspace v\defineq st,
\enspace \hbox{\stampatello both\enspace \sf and\enspace coprime}
$$
\par
\noindent
as ({\stampatello see the above})\hfill ${{\ell}\over s}\cdot r$ \sf, $st$ \sf by:\enspace $s|\ell$ \& $(t,\ell)=1$ $\Rightarrow $ $(t,s)=1$,\enspace  both $s,t$ \sf from $\mu(s)\mu(t)$\enspace {\stampatello twist};
\par
\noindent
{\stampatello and}\enspace $\left({{\ell}\over s}\cdot r,st\right)=1$\enspace {\stampatello from}\enspace $\left({{\ell}\over s},s\right)=1=\left({{\ell}\over s},t\right)$ (the former from $\ell$ \sf, the latter from $t$ coprime to $\ell$) {\stampatello and finally}\enspace $r\in )P(,\, st\in (P)$ $\Rightarrow $ $(r,st)=1$.
\par
\noindent
{\stampatello Thus}
$$
\hbox{\rm the\enspace map},\enspace \forall \ell\in \Pb\; \hbox{\stampatello fixed},
\qquad
(r,s,t)\in (3.2)\mapsto q=\ell rst^2\in \{ n\in \N : n\; \cf\} 
$$
\par
\noindent
(with $\in (3.2)$ meaning: variables like in $(3.2)$ up) {\stampatello is} an {\stampatello injective} function: \hfill {\stampatello we prove this, briefly}.
\par
\noindent
Our $q$ \lq \lq arrival\rq \rq \thinspace is\enspace $q=uv^2$ {\stampatello by the above}: {\stampatello if } $\exists (r_1,s_1,t_1),(r_2,s_2,t_2)\in (3.2)$,\hfill {\stampatello mapped to same } $uv^2$, 
\par
\noindent
{\stampatello then}:\enspace $u={{\ell}\over {s_1}}\cdot r_1={{\ell}\over {s_2}}\cdot r_2$ $\Rightarrow $ $s_2r_1=s_1r_2$, {\stampatello by $(P)$ \& $)P($ uniqueness of factorization} for natural numbers, {\stampatello which holds } $\forall P\in \Primes$ (fixed), we get \enspace $s_2=s_1$ \& $r_2=r_1$; 
{\stampatello and}:\enspace $v=s_1t_1=s_2t_2$ ({\stampatello with previous } $s_1=s_2$) $\Rightarrow $ $t_2=t_1$. Then, {\stampatello map is injective}.\hfill QED 
\par
We need only two ingredients, to calculate the $\lim_P$ of $(3.2)$, namely this {\stampatello map is one to one} and
$$
r\in )P(\enspace ,\enspace r>1
\enspace \Rightarrow \enspace  
r>P
\enspace \Rightarrow \enspace  
\ell rst^2 > P,
\enspace \forall \ell, r,s,t\in \N\enspace . 
\leqno{(r>P)}
$$
\par
\noindent
{\stampatello In fact}, $\forall \ell \in \Pb$  
$$
\Win_{\ell}^{(P)}\; \widetilde{F}=\sum_{s|\ell}\sum_{{t\in (P)}\atop {(t,\ell)=1}}\mu(s)\mu(t)G(\ell st^2)
+\sum_{s|\ell}\sum_{{t\in (P)}\atop {(t,\ell)=1}}\mu(s)\mu(t)\sum_{{r\in )P(}\atop {r>1}}\mu(r)G(\ell rst^2),
$$
\par
\noindent
where $r>1$ sum in RHS is, by $(r>P)$\enspace \&\enspace \lq \lq map's one to one\rq \rq, together with absolute convergence of Ramanujan Coefficient $G$, 
$$
\ll \sumflat_{{r\in )P(}\atop {r>1}}\sumflat_{s|\ell}\sumflat_{{t\in (P)}\atop {(t,\ell)=1}}\left|G(\ell rst^2)\right| 
 \ll \sum_{q>P}|G(q)| \buildrel{P}\over{\longrightarrow} 0,
\enspace \hbox{\stampatello uniformly} \enspace \forall \ell \in \Pb.
$$
\par
\noindent
{\stampatello Actually, if we fix $\ell \in \N$ \sf}, once\enspace $P\ge \ell$, we have\enspace $\ell \in \Pb$; in all, 
$$
\lim_P \;\Win_{\ell}^{(P)}\; \widetilde{F}=\lim_P \;\sum_{s|\ell}\sum_{{t\in (P)}\atop {(t,\ell)=1}}\mu(s)\mu(t)G(\ell st^2),
\enspace \hbox{\stampatello uniformly} \enspace \forall \ell \in \N^{\flat}.
$$
\par
\noindent
In turn, {\stampatello this entails, again by absolute convergence of Ramanujan Coefficient $G$} over \cf, 
$$
\lim_P \;\Win_{\ell}^{(P)}\; \widetilde{F}=\sum_{s|\ell}\sum_{{t=1}\atop {(t,\ell)=1}}^{\infty}\mu(s)\mu(t)G(\ell st^2),
\enspace \hbox{\stampatello uniformly} \enspace \forall \ell \in \N^{\flat},
$$
\par
\noindent
proving the other Formula.\hfill $\square$ 

\vfill
\eject

\par				
\noindent
It's time to prove our next result.
\medskip
\par
In order to prove Corollary 3.2, we need following general result. 
\smallskip
\par
\noindent {\bf Lemma 4.4}. ({\stampatello Ramanujan Sums with \sf argument \& square modulus})
\par
\noindent
{\it Let } $a,v\in \N$ {\it with } $a$\enspace {\stampatello square-free}. {\it Then}
$$
c_{v^2}(a)=\1_{v|a}\cdot v\cdot \mu(v). 
$$
\par
\noindent{\bf Proof}. From H\"older Formula $(1.6)$, 
$$
c_{v^2}(a)=\varphi(v^2){{\mu(v^2/(a,v^2))}\over {\varphi(v^2/(a,v^2))}}
 \buildrel{\thinspace a\thinspace \sf}\over{=\!=\!=} \varphi(v^2){{\mu(v^2/(a,v))}\over {\varphi(v^2/(a,v))}}
  \neq 0
\enspace \Longrightarrow \enspace v|a, 
$$
\par
\noindent
because otherwise\enspace $(a,v)<v$\enspace $\Rightarrow$\enspace ${v\over {(a,v)}}>1$\enspace $\Rightarrow$\enspace $\exists p\in \Primes$ with $p\left|{v\over {(a,v)}}\right.$ \& $p|v$\enspace $\Rightarrow$\enspace $p^2\left|{{v^2}\over {(a,v)}}\right.$\enspace $\Rightarrow$\enspace $\mu(v^2/(a,v))=0$; whence
$$
c_{v^2}(a)=\1_{v|a}{{\varphi(v^2)}\over {\varphi(v)}}\mu(v) 
$$
\par
\noindent
and use : 
$$
{{\varphi(v^2)}\over {v^2}}=\prod_{p|v}\left(1-{1\over p}\right)
  ={{\varphi(v)}\over v},
\enspace \forall v\in \N, \enspace 
$$
\par
\noindent
from $(1.14)$, to conclude.\hfill $\square$
\medskip
\par
We apply the Explicit Formula for normalized \enspace $\widetilde{C_{f,g_Q}}(N,a)$,\enspace coming from Theorem 3.6. 
\smallskip
\par
\noindent{\bf Proof of Corollary 3.2}. We only need to prove the asymptotic formula, soon after $\EF$, for the $3-$hp.s\enspace $C_{f,g_Q}(N,a)$,\enspace for fixed $a\in \N$ \sf, under both $\A$, $\P0$. We use $\EF$, to start: its LHS, {\stampatello thanks to}\enspace $P_0$\enspace  definition, see page 13, is supported on the \sf $\ell \in \left(P_0\right)$, too. From Corollary 3.1, $A_{f,g_Q}(N,a)$ (page 10 definition), the Analytic Part of $C_{f,g_Q}(N,a)$, is the same as its IPPification, $\widetilde{A_{f,g_Q}}(N,a)$, that's the Analytic Part of \thinspace $\widetilde{C_{f,g_Q}}(N,a)$, from the coincidence of respective Wintner-Carmichael coefficients. Hence, 
$$
\forall a\in \N^{\flat},
\enspace 
\widetilde{C_{f,g_Q}}(N,a)-A_{f,g_Q}(N,a)=N\sum_{q\;\cf}G(q)c_q(a)-\sumflat_{q\in \left(P_0\right)}{{\widehat{g_Q}(q)}\over {\varphi(q)}}\sum_{p\le N}f(p)c_q(p)c_q(a)
$$
\par
\noindent
(RHS $q\in \left(P_0\right)-$sum is over \sf from Hypothesis $(1)$ about $\widehat{g_Q}$) and normalizing, together with $\EF$, 
$$
\forall a\in \N^{\flat},
\enspace 
\widetilde{F}(a)-A_{f,g_Q}(N,a)/N=\sum_{q\;\cf}G(q)c_q(a)-\sumflat_{q\in \left(P_0\right)}\sum_{s|q}\sum_{{t=1}\atop {(t,q)=1}}^{\infty}\mu(s)\mu(t)G(qst^2)c_q(a). 
$$
\par
\noindent
Hereafter we {\stampatello fix $a\in \N$ square-free}. Then, applying $(2.1)$ to RHS, first for $q$ (after, for $qst^2$, following), say, 
$$
\widetilde{F}(a)-A_{f,g_Q}(N,a)/N=\mathop{\sumflat_{u}\sumflat_{v}}_{(u,v)=1}G(uv^2)c_u(a)c_{v^2}(a)-\Sigma,
\leqno{(\ast)_{\to 0}}
$$
\par
\noindent
abbreviating, hereafter recall: $G$ is \cf supported, 
$$
\Sigma:=\sumflat_{q\in \left(P_0\right)}\sum_{s|q}\sum_{{t=1}\atop {(t,q)=1}}^{\infty}\mu(s)\mu(t)G(qst^2)c_q(a). 
$$
\par
\noindent
Equation's name: $(\ast)_{\to 0}$ is not by chance as we want to prove that above quantity $\to 0$, as $N\to \infty$, using also $\P0$, because we will also need that $P_0\to \infty$. First $(\ast)_{\to 0}$ RHS term is, applying Lemma 4.4 above, 
$$
\mathop{\sumflat_{u}\sumflat_{v}}_{(u,v)=1}G(uv^2)c_u(a)c_{v^2}(a)=\mathop{\sumflat_{u}\sumflat_{v|a}}_{(u,v)=1}v\mu(v)G(uv^2)c_u(a). 
$$

\vfill
\eject

\par				
\noindent
We'll use above equation soon. In the meanwhile, we treat $\Sigma$ above: the unique $(2.1)$ representation 
$$
qst^2={q\over s}\cdot (st)^2,
\enspace \hbox{\stampatello with}\enspace 
\mu^2\left({q\over s}\right)=1=\mu^2\left(st\right)
\enspace \hbox{\stampatello and}\enspace 
\left({q\over s},st\right)=1
$$
\par
\noindent
{\stampatello allows the change of variables}: $q/s=u$ \& $st=v$, whence, from $(P_0)$ divisor-closed, $s,u\in (P_0)$ :
$$
\Sigma \EqByDef \sumflat_{s\in (P_0)}\sumflat_{{u\in (P_0)}\atop {(u,s)=1}}c_s(a)c_u(a)\sumflat_{{v=1}\atop {{(v,u)=1}\atop {v\equiv 0(\!\!\bmod s)}}}^{\infty}\mu(s)\mu\left({v\over s}\right)G(uv^2) 
 = \sumflat_{u\in (P_0)}c_u(a)\sumflat_{{v=1}\atop {(v,u)=1}}^{\infty}\mu(v)G(uv^2)\sum_{{s\in (P_0)}\atop {s|v}}c_s(a), 
$$
\par
\noindent
because the {\stampatello other two constraints}, for $s|v$, i.e.: $s$ \sf {\stampatello and} $(s,u)=1$ {\stampatello follow from} $v$ \sf {\stampatello and} $(v,u)=1$. 
\par
\noindent
$\underline{\hbox{\stampatello If}}$ $v\in (P_0)$,\enspace $\underline{\hbox{\stampatello then}}$ \enspace${\displaystyle \sum_{{s\in (P_0)}\atop {s|v}}c_s(a)=\sum_{s|v}c_s(a)=\1_{v|a}\cdot v}$,\enspace $\forall a\in \N$,\enspace {\stampatello from} $(1.8)$; clearly, $s|v$ and $s\in (P_0)$, in general, are equivalent to\enspace $s\left|v_{(P_0)}\right.$ : $v_{(P_0)}$ {\stampatello is the $(P_0)-$component of}\thinspace $v$ (see page 27), whence other case. 
$\underline{\hbox{\stampatello If}}$ $v\notin (P_0)$,\enspace $\underline{\hbox{\stampatello then}}$ \enspace${\displaystyle \sum_{{s\in (P_0)}\atop {s|v}}c_s(a)=\1_{v_{(P_0)}|a}\cdot v_{(P_0)}}$, again by $(1.8)$; this sum is between $0$ and $a$, whence $O(a)$. Hence, 
$$
\Sigma=\sumflat_{u\in (P_0)}c_u(a)\left(\sumflat_{{v\in (P_0)}\atop {{(v,u)=1}\atop {v|a}}}v\mu(v)G(uv^2)+\sumflat_{{v\notin (P_0)}\atop {(v,u)=1}}\mu(v)G(uv^2)O(a)\right). 
$$
\par
\noindent
In all, above $(\ast)_{\to 0}$ becomes, since above $\Sigma$ first term has $uv\in (P_0)$ : 
$$
\widetilde{F}(a)-A_{f,g_Q}(N,a)/N=\mathop{\sumflat_{u}\sumflat_{v|a}}_{{(u,v)=1}\atop {uv\notin (P_0)}}v\mu(v)G(uv^2)c_u(a)
-\sumflat_{u\in (P_0)}c_u(a)\sumflat_{{v\notin (P_0)}\atop {(v,u)=1}}\mu(v)G(uv^2)O(a). 
$$
\par
\noindent
{\stampatello Bounding}: $|c_u(a)|\le (a,u)\le a$ (Lemma A.1 [CM] \& trivial bound), $c_u(a)\ll_a 1$, {\stampatello we get} LHS is 
$$
\ll_a \enspace 
\mathop{\sumflat_{u}\sumflat_{v}}_{{(u,v)=1}\atop {uv\notin (P_0)}}|G(uv^2)| + \sumflat_{u\in (P_0)}\sumflat_{{v\notin (P_0)}\atop {(v,u)=1}}|G(uv^2)| 
\enspace \ll_a \enspace 
\mathop{\sumflat_{u}\sumflat_{v}}_{{(u,v)=1}\atop {uv>P_0}}|G(uv^2)| + \sumflat_{u}\sumflat_{{v>P_0}\atop {(v,u)=1}}|G(uv^2)|, 
$$
\par
\noindent
in which, trivially, $n\notin (P_0)$\enspace $\Rightarrow $\enspace $n>P_0$ is applied once to $n=uv$ and once to $n=v$; this entails 
$$
\widetilde{F}(a)-A_{f,g_Q}(N,a)/N=
O_a\left(
\sum_{{q\;\cf}\atop {q>P_0}}|G(q)| + \sum_{{q\;\cf}\atop {q>P_0^2}}|G(q)|
\right)=o(1),
$$
\par
\noindent
as $N\to \infty$, from $\P0$ \& absolute convergence of Ramanujan coefficients in $\A$.\hfill $\square$
\medskip

\vfill
\eject

\par				
\noindent
In order to prove Corollary 3.3 we need the following difference-estimate. 
\smallskip
\par
\noindent {\bf Lemma 4.5}. ({\stampatello Contribution to Analytic Part has negligible $p^K-$term for $K>1$})
\par
\noindent
{\it Let } $a,N\in \N$ {\it with } {\stampatello large} $N$, {\it while } $a$\enspace {\stampatello is fixed}. {\it Then}
$$
A_{\Lambda,\Lambda_N}(N,a)-A_{\mu^2\,\Lambda,\Lambda_N}(N,a)=O_{\varepsilon}\left(N^{\varepsilon}a\sqrt{N}\right),
\quad
\hbox{\stampatello as } N\to \infty; 
$$
\par
\noindent
{\stampatello in particular}, {\it this} {\stampatello difference}, {\it once} {\stampatello fixed } $a\in \N$,\enspace {\stampatello is } $o(N)$. 
\par
\noindent{\bf Proof}. Since LHS counts $p^K$ contributions with $K>1$, 
$$
A_{\Lambda,\Lambda_N}(N,a)-A_{\mu^2\,\Lambda,\Lambda_N}(N,a)=\sum_{q\le N}{{\widehat{\Lambda_N}(q)}\over {\varphi(q)}}\sum_{1<K\le {{\log N}\over {\log 2}}}\sum_{p^K\le N}\left(\log p\right)c_q(p^K)c_q(a) 
$$
\par
\noindent
which is bounded, {\stampatello using} quoted Lemma in [CM], $|c_q(a)|\le (a,q)\le a$, $\forall a,q\in \N$ (compare  above), as 
$$
\left|A_{\Lambda,\Lambda_N}(N,a)-A_{\mu^2\,\Lambda,\Lambda_N}(N,a)\right|\ll_{\varepsilon} N^{\varepsilon}a \sum_{q\le N}{1\over {q^2}}\left(\sum_{p|q}\varphi(p)+\sqrt{N}\right),  
$$
\par
\noindent
{\stampatello from}: \enspace $\widehat{\Lambda_N}(q)\ll {{N^{\varepsilon}}\over q}$,\enspace $\forall q\le N$; \enspace ${1\over {\varphi(q)}}\ll {{\log \log q}\over q}$,\enspace $\forall q\ge 3$ [MoV,th.m 2.9]; use of $(1.6)$\enspace \hbox{\stampatello for both cases}: 
$$
p\not | q
\enspace \Rightarrow \enspace 
(p,q)=1
\enspace \Rightarrow \enspace 
c_q(p^K)=\mu(q),
\enspace \hbox{\stampatello while} \enspace
p|q
\enspace \Rightarrow \enspace 
c_q(p^K)=-\varphi(p)\mu(q),
\enspace \forall q\; \hbox{\stampatello square-free}. 
$$
\par
\noindent
(Of course, we have \lq \lq {\stampatello absorbed}\rq \rq, so-to-speak, {\stampatello $\log-$powers within} $N^{\varepsilon}$, for above bound.)
\smallskip
\par
\noindent
Since : \enspace ${\displaystyle \sum_{p|q}\varphi(p)\le \sum_{d|q}\varphi(d)=q }$,\enspace {\stampatello from } $(1.14)$, we get difference \quad ${\displaystyle O_{\varepsilon}\left(N^{\varepsilon}a\sqrt{N}\right), }$ \enspace {\stampatello as above}.\hfill $\square$ 
\medskip
\par
Our above assumption $\P0$ holds for the 3hp.s-Correlation $C_{\mu^2\,\Lambda,\Lambda_N}(N,a)$: we prove this, briefly.
\smallskip
\par
\noindent {\bf Lemma 4.6}. $\left(\hbox{\stampatello Wintner's Prime for the Correlation}\enspace C_{\mu^2\, \Lambda,\Lambda_N}(N,a) : \hbox{\stampatello  it satisfies $\P0$ above}\right)$
$$
P_0(C_{\mu^2\, \Lambda,\Lambda_N})=\max\{p\in \Primes : p\le N\}
 >N/2
\enspace \Rightarrow \enspace 
P_0(C_{\mu^2\, \Lambda,\Lambda_N})\to \infty,
\enspace \hbox{\stampatello as} \enspace 
N\to \infty.  
$$
\par
\noindent{\bf Proof}. We'll show only equation above. By Theorem 3.1, $C_{\mu^2\, \Lambda,\Lambda_N}$\enspace has Wintner Coefficients: 
$$
\Win_q C_{\mu^2\,\Lambda,\Lambda_N} = {{\widehat{\Lambda_N}(q)}\over {\varphi(q)}}\sum_{p\le N}(\log p)c_q(p), 
\qquad
\forall q\in \N, 
$$
\par
\noindent
that vanish $\forall q>N$, since this is true for $\widehat{\Lambda_N}$ by definition. Also, this is, in the dyadic range, for {\stampatello large} $N$, 
$$
{N\over 2}<q\le N
\enspace \Rightarrow \enspace 
\widehat{\Lambda_N}(q) \EqByDef -\,{{\mu(q)\log(q)}\over q},
\enspace \hbox{\stampatello that's} \enspace {{\log(q)}\over q}\neq 0,
\enspace \hbox{\stampatello for}\enspace {N\over 2}<q\le N\enspace \hbox{\stampatello prime}
$$
\par
\noindent
in which range, say, \lq \lq Bertand's Postulate\rq \rq \thinspace ensures in fact the existence of primes $q$ (even Hoisel's first result on Primes in Short Intervals, compare [MoV], proves this as a consequence). A moment's reflection on the definition (at page 13) of Wintner's Prime, for above Wintner coefficients, reduces to show: defining $P_0:=\max\{p\in \Primes : p\le N\}$, since $P_0\in (P_0)$, 
$$
\Win_{P_0} C_{\mu^2\,\Lambda,\Lambda_N}\neq 0.
$$
\par
\noindent
In turn, after proved above \enspace $\widehat{\Lambda_N}(P_0)\neq 0$, we are left with:
$$
\sum_{p\le N}(\log p)c_{P_0}(p)\neq 0.
\leqno{(4.5)}
$$
\par				
\noindent
This is a consequence of H\"older's $(1.6)$ for \enspace $c_{P_0}(p)$: 
$$
p\neq P_0\enspace \Rightarrow \enspace c_{P_0}(p)=\mu(P_0)=-1
\quad \hbox{\stampatello and}\quad 
p=P_0\enspace \Rightarrow \enspace c_{P_0}(p)=\varphi(P_0)=P_0-1, 
$$
\par
\noindent
whence
$$
\sum_{p\le N}(\log p)c_{P_0}(p)=-\sum_{{p\le N}\atop {p\neq P_0}}\log p+(P_0-1)\log P_0
 =P_0\log P_0-\sum_{p\le N}\log p
  \sim P_0\log P_0-N
   \sim P_0\log P_0,
$$
\par
\noindent
as $N\to \infty$, by a form for the $\PNT$ involving weight $\log p$ (see soon after Fact 1.6 in $\S1$), proving $(4.5)$.\hfill $\square$
\bigskip
\par
It's time to see why, with the help of previous Lemmata, Corollary 3.3 comes from Corollary 3.2. 
\smallskip
\par
\noindent{\bf Proof of Corollary 3.3}. We start from the {\stampatello definition} of our Correlation:
{\stampatello }{\stampatello }{\stampatello }{\stampatello }{\stampatello }{\stampatello }{\stampatello }{\stampatello }{\stampatello }{\stampatello }
$$
\forall a,N\in \N
\qquad
C_{\mu^{2}\, \Lambda,\Lambda_N}(N,a)\defineq \sum_{n\le N}\mu^2(n)\Lambda(n)\Lambda_N(n+a). 
$$
\par
\noindent
We leave the check: it satisfies $(0)$, $(1)$, $(2)$, the 3 hp.s, to the careful Reader. 
\par
By Lemma 4.6, it has $\P0$ above. Following [C4], esp. the {\stampatello Proof} of {\stampatello Th.m 5.3}, 
$$
\forall a,N\in \N, 
\quad
C_{\mu^2\, \Lambda,\Lambda_N}(N,a)= C_{\Lambda,\Lambda_N}(N,a)+O\left(\sqrt{N}L\log(N+a)\right),  
$$  
\par
\noindent
{\stampatello abbreviating } $L:=\log N$, here. {\stampatello From Corollary 3.2},
$$
\forall a\in \N \enspace \sf\enspace \hbox{\stampatello fixed},
\quad
C_{\mu^2\, \Lambda,\Lambda_N}(N,a)= A_{\mu^2\, \Lambda,\Lambda_N}(N,a)+o(N),  
\enspace
\hbox{\it as}
\enspace 
N\to \infty 
$$  
\par
\noindent
that, together with previous formula and Lemma 4.5, gives
$$
\forall a\in \N \enspace \sf\enspace \hbox{\stampatello fixed},
\quad
C_{\Lambda,\Lambda_N}(N,a)= A_{\Lambda,\Lambda_N}(N,a)+o(N),  
\enspace
\hbox{\it as}
\enspace 
N\to \infty. 
$$
\par
\noindent
In [C0] we started from the $\REEF$ for $C_{\Lambda,\Lambda_N}(N,a)$ : 
$$
C_{\Lambda,\Lambda_N}(N,a)=A_{\Lambda,\Lambda_N}(N,a),
\quad
(\hbox{\stampatello under}\enspace \DH, \hbox{\stampatello in}\enspace [{\rm C}0]; \hbox{\rm compare\enspace page\enspace }15)
$$
\par   
\noindent
so, {\stampatello actually} [C0] proves, {\stampatello unconditionally}, that (our $\hbox{\gothic S}(a)$ {\stampatello is $\hbox{\gothic S}_{\Lambda,\Lambda}(a)$ in } [C0]) : 
$$
\forall a\in \N\enspace \hbox{\stampatello fixed},
\quad
A_{\Lambda,\Lambda_N}(N,a)=\hbox{\gothic S}(a)N+O(Ne^{-c\sqrt{\log N}}),  
\enspace
\hbox{\it as}
\enspace 
N\to \infty. 
$$
\par
\noindent
Since \enspace $C_{\Lambda,\Lambda_N}(N,a)=C_{\Lambda,\Lambda}(N,a)+O(aL\log(N+a))$ \enspace {\stampatello from divisors' cut in} $\S2$, 
$$
\forall a\in \N \enspace \sf\enspace \hbox{\stampatello fixed},
\quad
C_{\Lambda,\Lambda}(N,a)=\hbox{\gothic S}(a)N+o(N),  
\enspace
\hbox{\it as}
\enspace 
N\to \infty. 
$$
\par
\noindent
(By {\stampatello partial summation} [T], this is exactly {\stampatello Conjecture B} [HL], {\stampatello for} \sf  $a=2k$.).\hfill $\square$

\vfill
\eject

\par				
\noindent{\bf 5. Odds \& Ends}
\bigskip
\par
\noindent
Going back to Pinch Lemma, next similar one will help us proving $(1.9)$ from $(1.7)$, as promised. 
\smallskip
\par
\noindent {\bf Lemma 5.1}. {\it Let } $d,m\in \N$. {\it Then, for } $\chi_0(\!\bmod\; m)$, 
$$
\sum_{t|d}\mu\left({d\over t}\right)\chi_0(t)=\sum_{{t|d}\atop {(t,m)=1}}\mu\left({d\over t}\right)
 =\1_{\kappa(d)|m}\mu(d). 
$$
\par
\noindent{\bf Remark 4.} See that above RHS$=\1_{d|m}\mu(d)$ : $d$ {\stampatello not square-free makes it vanish anyway}.\hfill $\diamond$
\smallskip
\par
\noindent{\bf Remark 5.} (Added in version3) Previous Lemma \& Remark help Th.m 3.6 Proof, too.\hfill $\diamond$
\smallskip
\par
\noindent{\bf Proof}. Write: $d=d_m\cdot d^*$, with $d_m:=\prod_{p|m}p^{v_p(d)}$ and $d^*:=\prod_{p\nondiv m}p^{v_p(d)}$, so $(d_m,d^*)=1$ and $d^*\in \Z_m^*$, giving \enspace $\sum_{t|d,(t,m)=1}\mu(d/t)=\mu(d_m)\sum_{t|d^*}\mu(d^*/t)=\mu(d_m)\1_{d^*=1}$, from $(1.2)$, and by definition $d^*=1\Rightarrow d_m=d$, so LHS above equals \enspace $\mu(d_m)\1_{d^*=1}=\mu(d)\1_{d^*=1}=\mu(d)\1_{p|d\Rightarrow p|m}=\mu(d)\1_{\kappa(d)|m}$, RHS above.\hfill $\square$    
\medskip
\par
We write $O(1)$, for a quantity bounded as $x\to \infty$, now, even if depending on other parameters, in our:
\medskip
\leftline{{\stampatello Proof} of : $(1.7)\enspace \Longrightarrow\enspace (1.9)$.} 
\medskip
$$
\sum_{a\le x}c_q(n+a)c_{\ell}(a)\buildrel{(1.7)}\over{=\!=\!=}\sum_{a\le x}\sum_{{d|q}\atop {d|n+a}}d\mu\left(q\over d\right)\sum_{{t|\ell}\atop {t|a}}t\mu\left({\ell\over t}\right)
 =\sum_{d|q}d\mu\left(q\over d\right)\sum_{t|\ell}t\mu\left({\ell\over t}\right)S_{d,t}(x), 
$$
\par
\noindent
where, setting\enspace $d':=d/(d,t)$, $t':=t/(d,t)$, $n':=n/(d,t)$ \enspace and \enspace $\overline{t'}(\bmod \;d')$ defined by $\overline{t'}\cdot t'\equiv 1(\bmod \;d')$, 
$$
S_{d,t}(x)\defineq \sum_{{a\le x}\atop {{a\equiv 0(\!\!\bmod t)}\atop {a\equiv -n(\!\!\bmod d)}}}1
 = \sum_{{b\le x/t}\atop {tb\equiv -n(\!\!\bmod d)}}1
  = \1_{(d,t)|n}\cdot \sum_{{b\le x/t}\atop {b\equiv -\overline{t'} \cdot n'(\!\!\bmod d')}}1
   = \1_{(d,t)|n}\cdot {{x(d,t)}\over {dt}}+O(1), 
$$
\par
\noindent
because : \enspace $\sum_{b\le B,\,b\equiv K(\!\!\bmod\,M)}1=B/M+O(1)$, whatever is \thinspace $K(\bmod\,M)$, gives, by Lemma 5.1 and $(1.2)$,  
$$
\lim_x {1\over x}\sum_{a\le x}c_q(n+a)c_{\ell}(a)=\sum_{d|q}d\mu\left(q\over d\right)\sum_{t|\ell}t\mu\left({\ell\over t}\right)\1_{(d,t)|n}\cdot {{(d,t)}\over {dt}}
 = \sum_{{v|(\ell,q)}\atop {v|n}}v\sum_{d|q}\mu\left(q\over d\right)\sum_{t|\ell,(d,t)=v}\mu\left({\ell\over t}\right) 
$$
$$
= \sum_{{v|(\ell,q)}\atop {v|n}}v\sum_{d'\left|{q\over v}\right.}\mu\left({q/v}\over {d'}\right)\sum_{t'\left|{\ell\over v}\right.,(t',d')=1}\mu\left({{\ell/v}\over {t'}}\right) 
 = \sum_{{v|(\ell,q)}\atop {v|n}}v\mu\left(\ell\over v\right)\sum_{{d'\left|{q\over v}\right.}\atop {d'\equiv 0(\!\!\bmod \ell/v)}}\mu\left({q/v}\over {d'}\right)
  = \sum_{{v|(\ell,q)}\atop {v|n}}v\mu\left(\ell\over v\right)\1_{q=\ell} , 
$$
\par
\noindent
which is, again by $(1.7)$, \enspace $\1_{q=\ell}c_{\ell}(n)$.\hfill $\square$ 
\par
We used above Remark in the Proof, now working $\forall \ell,q\in \N$. (We improved Lemma 5.1, w.r.t. previous version 1, having square-free $\ell,q$.) 

\bigskip

\par
Now, we add a \lq \lq plan\rq \rq, for future work, from a kind of view of present 3rd-ver., 3rd version, results. 
\medskip
\par
\noindent
We didn't write (yet...) the following, in [C2]; but, everything, in this 3rd-ver., {\stampatello started from: once fixed\enspace $P\in \Primes$,\enspace uniformly } $\forall F\in \C^{\N}$, (a combinatorial identity as it does not require any hypothesis...) 
$$
F(a)=\sum_{q\in \Pb}\left(\Win_q^{\Pb}\; F\right)c_q(a),
\quad \forall a\in \Pb
\leqno\PGrid
$$
\par
\noindent
{\stampatello writing} \enspace $\Win_q^{\Pb}\; F:=\Win_q^{(P)}\; \widetilde{F}\EqByDef \Win_q\; F_{\Pb}$ \enspace, which involves the \lq \lq {\stampatello someone we know} \rq \rq \enspace $F_{\Pb}$ \enspace, allowing to write above as 
$$
F_{\Pb}(a)=\sum_{q\in \Pb}\left(\Win_q\; F_{\Pb}\right)c_q\left(a_{\Pb}\right),
\quad \forall a\in \N
\leqno{(5.1)}
$$

\vfill
\eject

\par				
\noindent
a  \lq \lq{\stampatello trick}\rq \rq, say, we used in {\stampatello Th.m 3.6 Proof}. (However, and this is decisive, with coefficient $G(q)$ instead.) 
\par
Here, the {\stampatello $q-$th coefficient}
$$
\Win_q\; F_{\Pb}\EqByDef \sum_{{d\in \Pb}\atop {d\equiv 0(\!\!\bmod \, q)}}{{F'(d)}\over d}=0, 
\quad
\forall q\in \N\backslash \Pb
$$  
\par
\noindent
{\stampatello and} 
\par
\centerline{
${\displaystyle 
\Win_q\; F_{\Pb} }$ \enspace 
{\stampatello vanishes, eventually, but it's a fin.sum}:
\enspace 
$\forall q\in \Pb$. 
}  
\medskip
\par
\noindent
In all,the \lq \lq $\PGrid$\rq \rq, for\enspace $F_{\Pb}$,\enspace above $(5.1)$, is not-strictly-requiring $q\in \Pb$. 
\par
However, {\stampatello since $\Pb$ is a fin.set ($\forall P\in \Primes$ fixed), } $\PGrid$ \lq \lq {\stampatello sounds like}\rq \rq \thinspace {\stampatello a \REEF \enspace for } $F_{\Pb}$: the {\stampatello problem is} the {\stampatello argument: we $\underline{\hbox{\stampatello can't}}$ have $c_q(a)$ in place of $c_q\left(a_{\Pb}\right)$ as } $a\in \N$, {\stampatello now. Yes, } the {\stampatello solution is}: {\stampatello take}\enspace $a\in \Pb$,\enspace so 
$c_q(a)\EqByDef c_q\left(a_{\Pb}\right)$. {\stampatello In other words}, $(5.1)$ above is {\stampatello a } $\underline{\hbox{\stampatello local}}$ \REEF,
 for $F_{\Pb}$ : {\stampatello valid } $\forall a\in \Pb$.
\medskip
\par
{\stampatello All these ideas} form the {\stampatello background} for {\stampatello above results}, given in 3rd-ver.
\par
\noindent
They \lq \lq {\stampatello maturated}\rq \rq \thinspace from studying $P-${\stampatello Wintner coefficients}, see [C2], \& {\stampatello calculating them for $\widetilde{F}$, $\underline{\hbox{\stampatello not}}$ for $F$} (why? {\stampatello answer above}, p.24).
\par
Since, $\forall P\in \Primes$ {\stampatello fixed}, {\stampatello (P-Grid)s are $\Pb-$local (R.E.E.F.)s}, it is {\stampatello very natural to ask}: {\stampatello can we link} $\underline{\hbox{\stampatello local}}$ {\stampatello Wintner coefficients} $\underline{\hbox{\stampatello to}}$\enspace $\underline{\hbox{\stampatello global}}$ {\stampatello Ramanujan coefficients} ? Answer in Theorem 3.6 (at least, for $\widetilde{F}$ with {\stampatello absolutely conv. Ramanujan Coefficients}).
\par
Notice: in the question, we are not asking a link to global {\stampatello Wintner} coefficients; we don't, since, at least for $\BH-${\stampatello Correlations} (and hence {\stampatello for 3-hp.s ones}), we know that {\stampatello Wintner's Transform is fin.$^{\hbox{\stampatello ly}}$-supported} (a theme for the careful Reader, {\stampatello to think about}). 
\par
This is {\stampatello a kind of} summary, to have an eagle-flight {\stampatello view} of {\stampatello above} results, added for 3rd-ver.(yes, it took me \lq \lq a while\rq \rq, {\stampatello to realize} how to put pieces together...). We didn't mention ({\stampatello FAI}) equation [C2], containing the  {\stampatello Irregular Part} $I_{f,g_Q}(N,a)$ with [C2] definition in terms of our {\stampatello Irregular Series}; equivalently: 
$$
I_{f,g_Q}(N,a)\defineq A_{f,g_Q}(N,a)-C_{f,g_Q}(N,a)=o(N), 
\enspace \forall a\in \N^{\flat}\enspace \hbox{\stampatello fixed},\enspace \hbox{\stampatello as}\enspace N\to \infty,
\quad
\hbox{\stampatello in}\enspace \hbox{\stampatello Corollary}\enspace 3.2.
\leqno{\52}
$$

\bigskip

\par
\noindent
{\stampatello Next step} (and I don't know how much will it take!) : {\stampatello give a $G:\N\rightarrow \C$, abs.conv.Ramanujan Coefficient}, more or less {\stampatello explicit for Corollary $3.2$, }, {\stampatello realizing $\A$ assumption.Maybe,} $(5.2)$ {\stampatello is a good start}.

\bigskip

That's all, folks!

\vfill
\eject

\par				
\centerline{\stampatello Bibliography}

\bigskip

\item{[Ca]} R. D. Carmichael, {\sl Expansions of arithmetical functions in infinite series}, Proc. London Math. Society {\bf 34} (1932), 1--26.
\smallskip
\item{[C0]} G. Coppola, {\sl An elementary property of correlations}, Hardy-Ramanujan J. {\bf 41} (2018), 65--76.
\smallskip
\item{[C1]} G. Coppola, {\sl Recent results on Ramanujan expansions with applications to correlations}, Rend. Sem. Mat. Univ. Pol. Torino {\bf 78.1} (2020), 57--82. 
\smallskip
\item{[C2]} G. Coppola, {\sl On Ramanujan smooth expansions for a general arithmetic function}, arXiv:2407.19759v1 
\smallskip
\item{[C3]} G. Coppola, {\sl Good Ramanujan Expansions: a suitably enhanced decay of coefficients has important consequences}, arXiv:2509.24456v3 (Third Version) 
\smallskip
\item{[CG]} G. Coppola and L. Ghidelli, {\sl Convergence of Ramanujan expansions, I [Multiplicativity on Ramanujan clouds]}, arXiv:2009.14121v1 
\smallskip
\item{[CM]} G. Coppola and M. Ram Murty, {\sl Finite Ramanujan expansions and shifted convolution sums of arithmetical functions, II}, J. Number Theory {\bf 185} (2018), 16--47. 
\smallskip
\item{[D]} H. Davenport, {\sl Multiplicative Number Theory}, 3rd ed., GTM 74, Springer, New York, 2000. 
\smallskip
\item{[De1]} H. Delange, {\sl On Ramanujan expansions of certain arithmetical functions}, Acta Arith., {\bf 31} (1976), 259--270.
\smallskip
\item{[De2]} H. Delange, {\sl On a formula for almost-even arithmetical functions}, Illinois J. Math., {\bf 31} (1987), 24--35.
\smallskip
\item{[H]} G.H. Hardy, {\sl Note on Ramanujan's trigonometrical function $c_q(n)$ and certain series of arithmetical functions}, Proc. Cambridge Phil. Soc. {\bf 20} (1921), 263--271.
\smallskip
\item{[HL]} G.H. Hardy and J.E. Littlewood, {\sl SOME PROBLEMS OF 'PARTITIO NUMERORUM'; III: ON THE EXPRESSION OF A NUMBER AS A SUM OF PRIMES.} Acta Mathematica {\bf 44} (1923), 1--70. 
\smallskip
\item{[IKo]} H. Iwaniec and E. Kowalski, {\sl Analytic Number Theory}, American Mathematical Society Colloquium Publications, 53. American Mathematical Society,Providence,RI, 2004. xii+615pp.ISBN:0-8218-3633-1
\smallskip
\item{[MoV]} H.L. Montgomery and R. Vaughan, {\sl Multiplicative Number Theory. I. Classical Theory}, Cambridge Studies in Advanced Mathematics, {97}, Cambridge University Press, 2007. 
\smallskip
\item{[M]} M. Ram Murty, {\sl Ramanujan series for arithmetical functions}, Hardy-Ramanujan J. {\bf 36} (2013), 21--33. Available online 
\smallskip
\item{[R]} S. Ramanujan, {\sl On certain trigonometrical sums and their application to the theory of numbers}, Transactions Cambr. Phil. Soc. {\bf 22} (1918), 259--276.
\smallskip
\item{[ScSp]} W. Schwarz and J. Spilker, {\sl Arithmetical Functions}, Cambridge University Press, 1994.
\smallskip
\item{[T]} G. Tenenbaum, {\sl Introduction to Analytic and Probabilistic Number Theory}, Cambridge Studies in Advanced Mathematics, {46}, Cambridge University Press, 1995. 
\smallskip
\item{[W]} A. Wintner, {\sl Eratosthenian averages}, Waverly Press, Baltimore, MD, 1943. 

\bigskip
\bigskip
\bigskip

\par
\leftline{\tt Giovanni Coppola - Universit\`{a} degli Studi di Salerno (affiliation)}
\leftline{\tt Home address : Via Partenio 12 - 83100, Avellino (AV) - ITALY}
\leftline{\tt e-mail : giocop70@gmail.com}
\leftline{\tt e-page : www.giovannicoppola.name}
\leftline{\tt e-site : www.researchgate.net}

\bye